\documentclass{amsart}
\usepackage{hyperref}
\hypersetup{
    pdftitle={Supersymmetry of chiral de Rham complex II},    
    pdfauthor={Reimundo Heluani},     
    pdfkeywords={Calabi-Yau, chiral}, 
    colorlinks=true,       
    linkcolor=blue,          
    citecolor=blue,        
    filecolor=blue,      
    urlcolor=blue           
}

\usepackage{amssymb,amsmath,amsfonts,euscript,mathrsfs,amsthm}

\def\End{\operatorname{End}}
\def\cR{{\mathscr R}}
\def\vac{|0\rangle}
\def\cU{{\mathscr U}}
\def\fo{{\mathfrak o}}
\def\cH{{\mathscr H}}
\DeclareMathOperator{\Id}{Id}
\def\fl{{\mathfrak l}}
\def\fs{{\mathfrak s}}

\newtheorem{thm}{Theorem}[section]
\numberwithin{equation}{section}

\newtheorem{prop}[thm]{Proposition}
\newtheorem{lem}[thm]{Lemma}

\theoremstyle{definition}
\newtheorem{defn}[thm]{Definition}

\newtheorem{ex}[thm]{Example}

\theoremstyle{remark}

\newtheorem{rem}[thm]{Remark}

\begin{document}
\title{Supersymmetry of the chiral de Rham complex II:\\Commuting sectors.}
\author{Reimundo Heluani}
\address{Department of Mathematics, UC Berkeley\\ CA 94720 \\ USA}
\email{heluani@math.berkeley.edu}
\thanks{Supported by the Miller Institute for basic research in science.}
\date{}
\maketitle
\begin{center}
	To Victor G. Kac on the occasion of his 65th birthday.
\end{center}
\begin{abstract}
	We construct two commuting $N=2$ structures on the space of sections of the
	chiral de Rham complex (CDR) of a Calabi-Yau manifold. We use this extra
	supersymmetry to construct  a non-linear automorphism of CDR
	preserving these $N=2$ structures. Finally, we show how to extend these
	results to construct a commuting pair of $N=4$ super-vertex
	algebras in the Hyper-K\"ahler case.
\end{abstract}
\section{Introduction}
In \cite{malikov}, the authors introduced a sheaf of vertex superalgebras
$\Omega_M^{\mathrm{ch}}$ attached
to any smooth manifold $M$, called the \emph{chiral de Rham complex} of $M$. The
sheaf cohomology $H^*(M, \Omega_M^{\mathrm{ch}})$ of $\Omega_M^\mathrm{ch}$ is
related to the chiral algebra of the half-twisted $\sigma$-model with target $M$, a
quantum field theory associated to $M$ (see \cite{FL}, \cite{kapustin},
\cite{Witten} and more recently \cite{frenkelnekrasov}). It was shown in
\cite{malikov} that in the holomorphic setting, when $M$ carries a global
holomorphic volume form, this super vertex algebra carries $N=2$ supersymmetry.

In
\cite{heluani2} the authors gave a \emph{superfield} formulation of CDR and studied
its properties in the $C^\infty$ setting. It was
shown there
that one can associate explicit sections of $\Omega_M^\mathrm{ch}$ to given
geometric tensors on $M$. For example, to a Riemmanian
metric $g_{ij}$ one associates a superfield $H$ (cf. \cite[7.4.1]{heluani2}) that
generates $N=1$ supersymmetry (that is, $H$ is a Neveu-Schwarz superconformal
vector). Moreover, any complex structure on $M$ gives rise to another superfield
$J$. The main result in
\cite{heluani2} states that the superfields $H$ and $J$ of the super vertex algebra $H^*(M, \Omega_M^\mathrm{ch})$
generate $N=2$ supersymmetry if $M$ is
Calabi-Yau. Moreover, if $(M, g_{ij}, J_1,J_2,J_3)$ is a Hyper-K\"ahler manifold
with the three complex structures $J_i$, $i=1,\dots, 3$ satisfying the quaternion
relations, then the corresponding superfields $\{H, J_i\}$ generate an $N=4$ super
vertex algebra. 

In this article, we show that on a K\"ahler manifold $(M, g, \omega)$, one can associate another
superfield $J_\omega$ to the K\"ahler form. It turns out that $H$ and $J_\omega$
generate $N=2$ supersymmetry on any K\"ahler manifold (see Theorem
\ref{thm:anothern=2}). Moreover, when $M$ is Calabi-Yau, the super-fields $H$, $J$
and $J_\omega$ generate two commuting copies of $N=2$ (see Theorem \ref{thm:n=2,2}). 

Recall from \cite{malikov} that, at least locally, $\Omega_M^\mathrm{ch}$ is
isomorphic to a $bc-\beta\gamma$ system in $\dim_\mathbb{R} M$ generators. 
Specifically, to a coordinate system $\{x_i\}$ on $M$, one associates local
sections $\{b^i, a_i, \phi^i, \psi_i\}$, $i=1,\dots, \dim_{\mathbb{R}} M$ of
$\Omega_M^{\mathrm{ch}}$. The fields $b^i$ transform as the coordinates $x_i$ in
$M$ do. The fields $\phi^i$ transform as the differential forms $dx^i$ do. The
fields $\psi_i$ transform as the vector fields $\tfrac{\partial}{\partial x_i}$ do,
and finally, the fields $a_i$ transform in a complicated way in order to cancel the
anomalies in the OPEs of these fields \cite{malikov}. The main idea of \cite{heluani2} is to use the fact that
there exists a \emph{supersymmetry} on the $bc-\beta\gamma$ system generated by
\begin{equation*}
	b^i \mapsto \phi^i, \qquad \psi_i \mapsto a_i,
\end{equation*}
in order to construct superfields from these pair of super-partners (cf.
\cite{heluani3}). In the Calabi-Yau case, the existence of the extra $N=2$
supersymmetry generated by the K\"ahler form $\omega$ allows us to interchange the
roles of differential forms and vector fields in the above picture. Indeed, one can
use the metric $g_{ij}$ and its inverse $g^{ij}$ to identify the tangent and
cotangent bundles of $M$ and we obtain a new supersymmetry of the $bc-\beta\gamma$
system that is generated by
\begin{equation*}
	b^i \mapsto \psi^i := \sum_j g^{ij} \psi_j, \qquad \phi_i:= \sum_j g_{ij}
	\phi^j \mapsto \tilde{a}_i
\end{equation*}
where the expression for $\tilde{a}_i$ is complicated and involves explicitly the
Levi-Civita connection associated to $g_{ij}$ (see Theorem
\ref{thm:isomorphismstrange}). This allows us to construct an automorphism of
$\Omega_M^{\mathrm{ch}}$ as a sheaf of super-vertex algebras on a Calabi-Yau
manifold $M$. Moreover, this automorphism preserves the two commuting $N=2$
structures associated to the complex and K\"ahler structures of $M$ (see
Proposition \ref{prop:automorphimaction}). 

Finally, when $M$ is Hyper-K\"ahler, we can construct three sections $J_i$ of
$\Omega_M^\mathrm{ch}$ associated to the three complex structures on $M$, and also
we can construct three sections $J_{\omega_i}$ associated to the corresponding
K\"ahler structures. We prove that these superfields, together with $H$, generate
two commuting copies of the $N=4$ super-vertex algebra.

The organization of this article is as follows. In section \ref{classical} we
recollect some results on vertex superalgebras following \cite{kac:vertex}. In
section \ref{sec:2} we recall some results and notation on SUSY vertex algebras
from \cite{heluani3}. In section \ref{sec:courant} we
introduce the chiral de Rham complex of $M$ as a sheaf of SUSY vertex algebras. We
rederive in the superfield formulation of \cite{heluani3} some of the results of
\cite{bressler1} and \cite{gerbes2}. In particular, we show how in this case, the axioms of
the SUSY Lambda bracket of \cite{heluani3} correspond to axioms of a Courant
algebroid.  In section \ref{sec:anothern=2} we construct a (family of)  $N=2$ superconformal
structure of central charge $3 \dim_\mathbb{R} M$ on any K\"ahler manifold $M$ extending the
$N=1$ structure corresponding to the Riemmanian metric on $M$. We include in this
section a Lemma showing that on any symplectic manifold we can associate another
$N=2$ structure to the symplectic form.  In section \ref{sec:n=2,2}, we show that
the two different $N=2$ structures corresponding to the complex and K\"ahler
structures on $M$ generate two commuting copies of the $N=2$ super vertex algebra
of central charge $3 \dim_\mathbb{C} M$. In section \ref{sec:n=4,4} we show how to
extend some of these results to the Hyper-K\"ahler case, in particular, we show
that there is a pair of commuting $N=4$ structures of central charge $3
\dim_{\mathbb{C}} M$ on $M$. In the Appendix we include
the technical computations and proofs of the main Theorems.

For an introduction to SUSY vertex algebras in general we refer the reader to
\cite{heluani3}. For the SUSY vertex algebra approach to the chiral de Rham
complex, we refer the reader to \cite{heluani2}.
Unless otherwise noted, all vector spaces, vector bundles, etc. are assumed to be
over the complex numbers. We will assume sums over repeated indexes. Indexes
$\alpha, \beta, \dots$ will run over holomorphic coordinates, indexes
$\bar{\alpha}, \bar{\beta}, \dots$ will run over anti-holomorphic coordinates and
indexes $i,j, \dots$ will run over arbitrary coordinates. We will raise and lower
indexes with the metric tensor and its inverse.

\textbf{Acknowledgements:} the author would like to thank M. Gualtieri for
explaining the symmetries of the Courant bracket to him. He would like to thank E.
Frenkel for numerous explanations -- including the suggestion of complexifying the
K\"ahler form, M. Aldi and A. Linshaw for interesting
discussions and M. Szczesny for his generosity in sharing his insights.

\section{Vertex superalgebras} \label{classical}

In this section, we review the definition of vertex superalgebras,
as presented in \cite{kac:vertex}.
    Given a vector space $V$, an \emph{$\End(V)$-valued field} is a formal
    distribution of the form
    \begin{equation*}
        A(z) = \sum_{n \in \mathbb{Z}} z^{-1-n} A_{(n)},\qquad A_{(n)} \in
        \End(V),
    \end{equation*}
    such that for every $v \in V$, we have $A_{(n)}v = 0$ for large enough $n$.

    \begin{defn}
    A vertex super-algebra consists of the data of a super vector space $V$,
    an even vector $\vac \in V$ (the vacuum vector),
     an even endomorphism $T $, and a parity preserving linear map $A \mapsto Y(A,z)$ from
     $V$ to $\End(V)$-valued fields (the state-field correspondence). This
     data should satisfy the following set of axioms:
    \begin{itemize}
    \item Vacuum axioms:
        \begin{equation*}
            \begin{aligned}
            Y(\vac, z) &= \Id \\
            Y(A, z) \vac &= A + O(z) \\
            T \vac &= 0
            \end{aligned}
        \end{equation*}
    \item Translation invariance:
        \begin{equation*}
            \begin{aligned}
                {[}T, Y(A,z)] &= \partial_z Y(A,z)
        \end{aligned}
        \end{equation*}
    \item Locality:
        \begin{equation*}
            (z-w)^n [Y(A,z), Y(B,w)] = 0 \qquad n \gg 0
        \end{equation*}
     \end{itemize}
(The notation $O(z)$ denotes a power series in $z$ without constant
term.)
    \label{defn:1.3}
\end{defn}
    Given a vertex super-algebra $V$ and a vector $A \in V$, we expand the fields
    \begin{equation*}
        Y(A,z) = \sum_{{j \in \mathbb{Z}}} z^{-1-j}
        A_{(j)}
    \end{equation*}
    and we call the endomorphisms $A_{(j)}$ the \emph{Fourier modes} of
    $Y(A,z)$. Define now the operations:
    \begin{equation*}
        \begin{aligned}
            {[}A_\lambda B] &= \sum_{{j \geq 0}}
            \frac{\lambda^{j}}{j!} A_{(j)}B \\
            A B &= A_{(-1)}B
        \end{aligned}
    \end{equation*}
    The first operation is called the $\lambda$-bracket and the second is
    called the \emph{normally ordered product}.
     The $\lambda$-bracket contains all of the information about the commutators between the Fourier coefficients of fields in $V$. 

\subsection{The $N=1$, $N=2$, and $N=4$ superconformal vertex algebras}

In this section we review the standard description of the $N=1,2,4$
superconformal vertex algebras. In section \ref{sec:2}, the same
algebras will be described in the SUSY vertex algebra formalism.

\begin{ex}{\bf The $N=1$ (Neveu-Schwarz)
    superconformal vertex algebra} \label{N1ex}

    The $N=1$ superconformal vertex algebra (\cite{kac:vertex}) of central charge $c$ is generated by two fields: $L(z)$, an even field of conformal weight $2$, and $G(z)$, an odd primary field of conformal weight $\frac{3}{2}$, with the $\lambda$-brackets
\begin{equation*}
{[L}_\lambda L] = (T + 2\lambda) L + \frac{c \lambda^3}{12}
    \label{eq:1}
\end{equation*}
\[
{[L}_\lambda G] = (T+\frac{3}{2} \lambda) G
\]
\[
{[G}_\lambda G]  = 2L + \frac{c \lambda^2}{3}
\]
$L(z)$ is called the Virasoro field.
\end{ex}

\begin{ex}{\bf The $N=2$ superconformal vertex algebra} \label{N2ex}

The $N=2$ superconformal vertex algebra of central charge $c$ is
generated by the Virasoro field $L(z)$ with $\lambda$-bracket
(\ref{eq:1}), an even primary field $J(z)$ of conformal weight $1$,
and two odd primary fields $G^{\pm}(z)$ of conformal weight
$\frac{3}{2}$, with the $\lambda$-brackets \cite{kac:vertex}
\begin{equation*}
    {[}L_\lambda J] = (T + \lambda) J
\end{equation*}
\begin{equation*}
    [L_\lambda G^\pm] = \left(  T + \frac{3}{2}\lambda \right) G^\pm
\end{equation*}
\begin{xalignat*}{2}
    {[J}_\lambda G^\pm] &= \pm G^\pm & [J_\lambda J] &= \frac{c}{3}\lambda \\
    {[G^+}_\lambda G^-] &= L + \frac{1}{2} TJ + \lambda J + \frac{c}{6}\lambda^2
    & {[G^\pm}_\lambda G^\pm] &= 0
\end{xalignat*}

\end{ex}

\begin{ex}{\bf The ``small'' $N=4$ superconformal vertex algebra} \label{N4ex}

    The even part of this vertex algebra is generated by the Virasoro field
    $L(z)$ and three primary fields of conformal weights $1$, $J^0$, $J^+$ and $J^-$. The odd
    part is generated by four primary fields of conformal weight $3/2$, $G^\pm$
    and $\bar{G}^\pm$. The remaining
(non-vanishing)
$\lambda$-brackets are (cf \cite[page 36]{kacwakimoto1})
\begin{xalignat*}{2}
    {[J^0}_\lambda J^\pm] &= \pm 2 J^\pm & {[J^0}_\lambda J^0] &=
    \frac{c}{3}\lambda \\
    {[J^+}_\lambda J^-] &= J^0 + \frac{c}{6} \lambda & {[J^0}_\lambda G^\pm] &= \pm
    G^\pm \\
    {[J^0}_\lambda \bar{G}^\pm] &= \pm \bar{G}^\pm & {[J^+}_\lambda G^-] &= G^+
    \\ {[J^-}_\lambda G^+]&= G^- & {[J^+}_\lambda \bar{G}^-] &= - \bar{G}^+ \\
    {[J^-}_\lambda \bar{G}^+] &= - \bar{G}^- & {[G^\pm}_\lambda \bar{G}^\pm]&=
    (T + 2 \lambda) J^\pm \\
    {[G^\pm}_\lambda \bar{G}^{\mp}]&= L \pm \frac{1}{2}T J^0 \pm \lambda J^0 +
    \frac{c}{6} \lambda^2
\end{xalignat*}
\end{ex}
(Note that the $J$ currents form an $\fs\fl_2$ current
algebra.)

\section{SUSY vertex algebras}\label{sec:2} In this section we collect some
results on SUSY vertex algebras from \cite{heluani3}. 

\subsection{Structure theory of SUSY VAs}
    Introduce formal variables $Z=(z,\theta)$ and $W =
    (w,\zeta)$, where $\theta, \zeta$ are odd
    anti-commuting variables and $z, w$ are even commuting variables.
    Given an integer $j$ and $J = 0$ or $1$ we put $Z^{j|J} = z^j \theta^J$.

    Let $\cH$ be the superalgebra generated by $\chi, \lambda$ with the relations
    $[\chi, \chi] = - 2 \lambda$, where $\chi$ is
    odd and $\lambda$ is even and central. We will consider another set of
    generators $-S, -T$ for $\cH$ where $S$ is odd, $T$ is central, and $[S, S]
    = 2 T$. Denote $\Lambda = (\lambda, \chi)$,
    $\nabla = (T, S)$, $\Lambda^{j|J} = \lambda^j \chi^J$ and $\nabla^{j|J} =
    T^j S^J$.

    Given a super vector space $V$ and a vector $a \in V$, we will denote by
    $(-1)^a$ its parity.
    Let $U$ be a vector space, a $U$-valued formal distribution is an
    expression of the form
    \begin{equation*}
        \sum_{\stackrel{j \in \mathbb{Z}}{J = 0,1}} Z^{-1-j|1-J} w_{(j|J)}
        \qquad w_{(j|J)} \in U.
    \end{equation*}
    The space of such distributions will be denoted by $U[ [Z, Z^{-1}] ]$. If
    $U$ is a Lie algebra we will say that two such distributions $a(Z), \,
    b(W)$ are
    \emph{local} if
    \begin{equation*}
        (z - w)^n [a(Z), b(W)] = 0 \qquad n \gg 0.
    \end{equation*}
    The space of distributions such that only finitely many negative powers
    of $z$ appear (i.e. $w_{(j|J)} = 0$ for large enough $j$) will be denoted
    $U( (Z ))$. In the case when $U = \End(V)$ for another vector space $V$,
    we will say that a distribution $a(Z)$ is a \emph{field} if $a(Z)v \in V(
    (Z ))$ for all $v \in V$.
    \begin{defn}[\cite{heluani3}]
    An $N_K=1$ SUSY vertex algebra consists of the data of a vector space $V$,
    an even vector $\vac \in V$ (the vacuum vector), an odd endomorphism
    $S$ (whose square is an even endomorphism we denote $T$),
    and a parity preserving linear map $A \mapsto Y(A,Z)$ from
     $V$ to $\End(V)$-valued fields (the state-field correspondence). This
     data should satisfy the following set of axioms:
    \begin{itemize}
    \item Vacuum axioms:
        \begin{equation*}
            \begin{aligned}
            Y(\vac, Z) &= \Id \\
            Y(A, Z) \vac &= A + O(Z) \\
            S \vac &= 0
            \end{aligned}
        \end{equation*}
    \item Translation invariance:
        \begin{equation*}
            \begin{aligned}
            {[} S, Y(A,Z)] &= (\partial_\theta - \theta \partial_z)
            Y(A,Z)\\
            {[}T, Y(A,Z)] &= \partial_z Y(A,Z)
        \end{aligned}
        \end{equation*}
    \item Locality:
        \begin{equation*}
            (z-w)^n [Y(A,Z), Y(B,W)] = 0 \qquad n \gg 0
        \end{equation*}
     \end{itemize}
    \label{defn:2.3}
\end{defn}
\begin{rem}
    Given the vacuum axiom for a SUSY vertex algebra, we will use the state
    field correspondence to identify a vector $A \in V$ with its corresponding
    field $Y(A,Z)$.
    \label{rem:nosenosenose}
\end{rem}
    Given a $N_K=1$ SUSY vertex algebra $V$ and a vector $A \in V$, we expand the fields
    \begin{equation*}
        Y(A,Z) = \sum_{\stackrel{j \in \mathbb{Z}}{J = 0,1}} Z^{-1-j|1-J}
        A_{(j|J)}
    \end{equation*}
    and we call the endomorphisms $A_{(j|J)}$ the \emph{Fourier modes} of
    $Y(A,Z)$. Define now the operations:
    \begin{equation}
        \begin{aligned}
            {[}A_\Lambda B] &= \sum_{\stackrel{j \geq 0}{J = 0,1}}
            \frac{\Lambda^{j|J}}{j!} A_{(j|J)}B \\
            A B &= A_{(-1|1)}B
        \end{aligned}
        \label{eq:2.4.2}
    \end{equation}
    The first operation is called the $\Lambda$-bracket and the second is
    called the \emph{normally ordered product}.
\begin{rem}
    As in the standard setting, given a SUSY VA $V$ and a vector $A \in V$, we
    have:
    \begin{equation*}
        Y(TA, Z) = \partial_z Y(A,Z) = [T, Y(A,Z)]
    \end{equation*}
    On the other hand, the action of the derivation $S$ is described
    by:
    \begin{equation*}
        Y(SA,Z) = \left( \partial_\theta + \theta \partial_z \right) Y(A,Z)
        \neq [S, Y(A,Z)].
    \end{equation*}
    \label{rem:caca5}
\end{rem}
The relation with the standard field formalism is as follows.
Suppose that $V$ is a vertex superalgebra  as defined in section
\ref{classical}, together with a homomorphism from the $N=1$
superconformal vertex algebra in example \ref{N1ex}. $V$ therefore
possesses an even vector $\nu$ of conformal weight $2$, and an odd
vector $\tau$ of conformal weight $\frac{3}{2}$, whose associated
fields
\begin{equation*}
    \begin{aligned}
 Y(\nu,z) &= L(z) = \sum_{n \in \mathbb{Z}} L_n z^{-n-2} \\
Y(\tau,z) &= G(z) = \sum_{n \in 1/2 + \mathbb{Z}} G_n z^{-n - \frac{3}{2}}
\end{aligned}
\end{equation*}
have the $\lambda$-brackets as in example \ref{N1ex}, and where we
require $G_{-1/2}=S$ and $L_{-1}=T$. We can then endow $V$ with the
structure of an $N_K=1$ SUSY vertex algebra via the state-field
correspondence \cite{kac:vertex}
\[
Y(A,Z) = Y^{c} (A,z) + \theta Y^{c}(G_{-1/2} A, z)
\]
where we have written $Y^{c}$ to emphasize that this is the
``classical" state-field (rather than state--superfield)
correspondence in the sense of section \ref{classical}.

(Note however that there exist $SUSY$ vertex algebras without such a map
from the $N=1$ superconformal vertex algebra.)
\begin{defn}
    Let $\cH$ be as before. \emph{An $N_K=1$ SUSY Lie
    conformal  algebra} is a $\cH$-module $\cR$ with an operation
    $[\,_{\Lambda}\,]: \cR \otimes \cR \rightarrow \cH
    \otimes \cR$ of degree
    $1$ satisfying:
    \begin{enumerate}
        \item Sesquilinearity
            \begin{equation*}
                [S a_\Lambda b] =  \chi [a_\Lambda b]
                \qquad [a_\Lambda S b] = -(-1)^{a} \left(S
                + \chi
                \right) [a_\Lambda b]
            \end{equation*}
        \item Skew-Symmetry:
            \begin{equation*}
                [b_\Lambda a] =  (-1)^{a b} [b_{-\Lambda -
                \nabla} a]
            \end{equation*}
            Here the bracket on the right hand side is computed as
            follows: first compute $[b_{\Gamma}a]$, where $\Gamma =
            (\gamma, \eta)$ are generators of $\cH$ super commuting
            with $\Lambda$, then replace $\Gamma$ by $(-\lambda - T,
            -\chi - S)$.
        \item Jacobi identity:
            \begin{equation*}
                [a_\Lambda [b_\Gamma c]] = -(-1)^{a} \left[
                [ a_\Lambda b]_{\Gamma + \Lambda} c \right] +
                (-1)^{(a+1)(b+1)} [b_\Gamma [a_\Lambda c]]
            \end{equation*}
            where the first bracket on the right hand side is computed as in Skew-Symmetry
            and the identity is an identity in $\cH^{\otimes 2} \otimes\cR$.
    \end{enumerate}
    \label{defn:k.conformal.1}
    \end{defn}
    Given an $N_K=1$ SUSY VA, it is canonically an $N_K=1$ SUSY Lie conformal algebra with
    the bracket defined in (\ref{eq:2.4.2}). Moreover, given an $N_K=1$ Lie
    conformal algebra $\cR$, there exists a unique $N_K=1$ SUSY VA called the
    \emph{universal enveloping SUSY vertex algebra of $\cR$} with the property
    that if $W$ is another $N_K=1$ SUSY VA and $\varphi : \cR \rightarrow W$ is a
    morphism of Lie conformal algebras, then $\varphi$ extends uniquely to a
    morphism $\varphi: V \rightarrow W$ of SUSY VAs.
    The operations (\ref{eq:2.4.2}) satisfy:
    \begin{itemize}
        \item Quasi-Commutativity:
            \begin{equation*}
                ab - (-1)^{ab} ba = \int_{-\nabla}^0 [a_\Lambda
                b] d\Lambda
            \end{equation*}
        \item Quasi-Associativity
            \begin{equation*}
                (ab)c - a(bc) = \sum_{j \geq 0}
                a_{(-j-2|1)}b_{(j|1)}c + (-1)^{ab} \sum_{j \geq 0}
                b_{(-j-2|1)} a_{(j|1)}c
            \end{equation*}
        \item Quasi-Leibniz (non-commutative Wick formula)
            \begin{equation*}
                [a_\Lambda bc ] = [a_\Lambda b] c + (-1)^{(a+1)b}b
                [a_\Lambda c] + \int_0^\Lambda [ [a_\Lambda
                b]_\Gamma c] d \Gamma
            \end{equation*}
    \end{itemize}
    where the integral $\int d\Lambda$ is $\partial_\chi \int d\lambda$. In
    addition, the vacuum vector is a unit for the normally ordered product
    and the endomorphisms $S, T$ are odd and even derivations respectively of
    both operations.
\subsection{Examples}
\begin{ex}
    Let $\cR$ be the free $\cH$-module generated by an odd vector $H$.
    Consider the following Lie conformal algebra structure in $\cR$:
    \begin{equation*}
        {[}H_\Lambda H] = (2T + \chi S + 3 \lambda) H
    \end{equation*}
    This is the \emph{Neveu-Schwarz} algebra (of central charge 0). This
    algebra admits a central extension of the form:
    \begin{equation*}
        {[}H_\Lambda H] = (2T + \chi S + 3\lambda) H + \frac{c}{3} \chi \lambda^2
    \end{equation*}
    where $c$ is any complex number. The associated universal enveloping
    SUSY VA is the \emph{Neveu-Schwarz} algebra of central charge
    $c$\footnote{Properly speaking, we consider the universal enveloping
    SUSY vertex algebra of $\cR \oplus \mathbb{C} C$ with $C$ central and $TC
    = S C = 0$ and then we quotient by the ideal generated by $C = c$ for
    any complex number $c$}. If we
    decompose the corresponding field
    \begin{equation*}
        H(z,\theta) = G(z) + 2 \theta L(z)
    \end{equation*}
    then the fields $G(z)$ and $L(z)$ satisfy the commutation relations of
    the well known $N=1$ super vertex algebra in example \ref{N1ex}.
    \label{ex:2.9}
\end{ex}
\begin{ex}
    Consider now the free $\cH$ module generated by even vectors $\{B^i\}_{i
    = 1}^n$ and odd vectors $\{\Psi_i\}_{i=1}^n$ where the only non-trivial
    commutation relations are:
    \begin{equation*}
	    {[B^i}_\Lambda \Psi_j] = \delta^{i}_{j} = {[\Psi_j}_\Lambda B^i]
    \end{equation*}
    Expand the corresponding fields as:
    \begin{equation*}
        B^i(z,\theta) = b^i(z) + \theta \phi^i(z) \qquad \Psi_i(z,\theta)
        = \psi_i(z) + \theta a_i(z)
    \end{equation*}
    then the fields $b^i$, $a_i$, $\phi^i$ and $\psi_i$ generate the
    $bc-\beta\gamma$ system as in \cite{malikov}.
    \label{ex:2.10}
\end{ex}
\begin{ex} \label{ex:2.11.a}
    The $N=2$ superconformal vertex algebra is generated by $4$
    fields \cite{kac:vertex}. In this context it is generated by two superfields -- an $N=1$ vector $H$ as in
    \ref{ex:2.9} and an even current $J$, primary of conformal weight $1$,
    that is:
    \begin{equation*}
        {[}H_\Lambda J] = (2T + 2\lambda + \chi S) J.
    \end{equation*}
    The remaining commutation relation is
    \begin{equation*}
        [J_\Lambda J] = - (H + \frac{c}{3} \lambda \chi).
    \end{equation*}
    Note that given
    the \emph{current} $J$ we can recover the $N=1$ vector $H$.
     In terms of the fields of Example \ref{N2ex}, $H, J$ decompose as

\begin{equation*}
    \begin{aligned}
         J(z,\theta) &= - \sqrt{-1}J(z) - \sqrt{-1}\theta \left( G^-(z) - G^+(z)
        \right) \\
        H(z,\theta) &= \left( G^+(z) + G^-(z) \right) + 2 \theta L(z)
    \end{aligned}
\end{equation*}

\end{ex}
\begin{ex} \label{ex:2.11}
    The ``small'' $N=4$ superconformal vertex algebra is a vertex algebra generated by 8
    fields \cite{kac:vertex}. In this formalism, it is generated by four
    superfields $H, J^i$, $i = 0,1,2$, such that each pair $(H, J^i)$ forms an
    $N=2$ SUSY VA as in the previous example and the remaining commutation
    relations are:
    \begin{equation*}
            {[}J^i_\Lambda J^j] =
            \varepsilon^{ijk} (S + 2\chi)J^k \qquad i \neq j
    \end{equation*}
    where $\varepsilon$ is the totally antisymmetric
    tensor.      In terms of the fields of Example \ref{N4ex}, $H, J^{i}$ decompose as

\begin{equation*}
    \begin{aligned}
        J^0(z,\theta) &= - \sqrt{-1} J^0(z) - \sqrt{-1} \theta \left(
        \bar{G}^-(z) - G^+(z)
        \right) \\
        J^1(z,\theta) &= \sqrt{-1} \left( J^+(z) + J^-(z) \right) +
        \sqrt{-1} \left( \bar{G}^+(z) - G^-(z) \right) \\
        J^2(z,\theta) &= \left(J^+(z) - J^-(z)\right) + \theta \left(
        \bar{G}^+(z) + G^-(z)
        \right) \\
        H(z,\theta) &= \left(G^+(z) + \bar{G}^-(z) \right) + 2 \theta L(z)
    \end{aligned}
\end{equation*}

\end{ex}

\section{Courant algebroids and the SUSY Lambda bracket} \label{sec:courant}

In this section we re-derive some of the results of \cite{bressler1} and
\cite{gerbes2} in the
language of SUSY vertex algebras. For a general introduction to Courant algebroids,
we refer the reader to
\cite{gualtieri1} and references therein. 



Let $M$ be a smooth manifold and denote by $T$ the complexified tangent bundle of
$M$.

\begin{defn}
	A Courant algebroid is a vector bundle $E$ over $M$, equipped with a
	nondegenerate symmetric bilinear form $\langle ,\rangle $ as
	well as a skew-symmetric bracket $[,]$ on
	$C^\infty(E)$ and with a smooth bundle map $\pi: E
	\rightarrow  T$ called the anchor. This induces a natural
	 differential operator $\mathcal{D}: C^\infty(M) \rightarrow
	 C^\infty(E)$ as $\langle \mathcal{D}f, A\rangle  = \tfrac{1}{2} \pi(A) f$ for
	all $f \in C^\infty (M)$ and $A \in
	C^\infty(E)$. 
 	 These structures should satisfy:
	\begin{enumerate}
		\item $\pi([A,B]) = [\pi(A), \pi(B)], \quad \forall A, B
			\in C^\infty(E)$. 
		\item The bracket $[,]$ should satisfy the following
			analog of the Jacobi identity. If we define the
			\emph{Jacobiator} as $\mathrm{Jac}(A,B,C) =
			[ [A,B], C] + [ [B, C], A ] + [ [ C,A],B]$.  And
			the \emph{Nijenhuis} operator 
			\[ \mathrm{Nij} (A,B,C) = \frac{1}{3} \left( \langle 
			[A,B], C\rangle  + \langle  [B,C], A\rangle  + \langle  [C, A], B\rangle 
			\right). \]
			Then the following must be satisfied:
			\[\mathrm{Jac} (A,B, C) = \mathcal{D} \left(
			\mathrm{Nij} (A, B, C)
			\right), \quad \forall A, B, C \in
			C^\infty(E) \]
		\item $[A, f B] = (\pi(A) f) B + f [A, B] - \langle A, B\rangle 
			\mathcal{D}
			f$, for all $A, B \in C^\infty(E)$ and
			$f \in C^\infty(M)$, 
		\item $\pi \circ \mathcal{D} = 0$, i.e. $\langle \mathcal{D} f,
			\mathcal{D} g\rangle  = 0,
			\quad \forall f,g \in C^\infty(M)$. 
		\item $\pi(A)\langle B, C\rangle  = \langle [A, B] + \mathcal{D}\langle A,B\rangle , C\rangle  + \langle B, [A, C]
			+ \mathcal{D} \langle A,C\rangle \rangle , \quad \forall A, B, C \in
			C^\infty(E)$. 
	\end{enumerate}
	\label{defn:1}
\end{defn}

\begin{ex}
	$E=T \oplus T^*$, $\langle,\rangle$ and $[,]$ are respectively the natural symmetric pairing and
	the Courant bracket defined as:
	\begin{equation*}
		\begin{aligned} \langle X + \zeta, Y + \eta\rangle &= \frac{1}{2} \left( i_X \eta + i_Y\zeta
\right).  \\
 {[X + \zeta}, Y + \eta] &= [X, Y] + \mathrm{Lie}_X \eta -
\mathrm{Lie}_Y \zeta - \frac{1}{2} d(i_X \eta - i_Y \zeta).
\end{aligned}
\end{equation*}
\label{ex:courant1}
\end{ex}

Let $(E, \langle ,\rangle , [,], \pi)$ be a Courant algebroid. 
Let  $\Pi E$ be the corresponding purely odd super vector bundle. we will abuse
notation and denote by $\langle, \rangle$ the corresponding
super-skew-symmetric bilinear form, and by $[,]$ the corresponding super-skew-symmetric
degree $1$ bracket on $\Pi E$. Similarly, we obtain an odd
differential operator $\mathcal{D}: C^\infty(M) \rightarrow 
C^\infty(\Pi E)$. If no confusion should arise, when $v$ is an element of a vector
space $V$, we will denote by the same symbol $v$ the corresponding element of $\Pi
V$, where $\Pi$ is the \emph{parity change operator}. 

Recall that $\cH$ is the
associative superalgebra generated by $S, T$ with relations $S^2 = T$
\cite{heluani3}.
Let $\cR$ be the $\cH$ module generated by $C^\infty(M) \oplus
 C^\infty(\Pi E)$ with the relation $S f = \mathcal{D} f$ for all $f \in
C^\infty(M)$. 

\begin{prop} The following endows $\cR$ with the structure of
an $N_K=1$ SUSY Lie conformal algebra:
\begin{equation}
	\begin{aligned}
		{[f}_\Lambda g] &= 0, \qquad \forall f,g \in
		C^\infty(M), \\
		{[A}_\Lambda B ] &= [A,B] + (2 \chi + \mathcal{D}) \langle A,
		B\rangle , \quad \forall A, B \in \Pi C^\infty(E), \\
		{[A}_\Lambda f] &= \pi(A) f, \qquad \forall A \in
		\Pi C^\infty(E), f \in C^\infty(M).
	\end{aligned}
	\label{eq:lambdadef}
\end{equation}
\label{prop:susycommutators}
\end{prop}

\begin{proof}
	Let us check first that the Lambda bracket is well defined:
\begin{equation*}
	[A_\Lambda \mathcal{D} f] = [A, \mathcal{D} f]
	+ (2 \chi + \mathcal{D}) \langle A, \mathcal{D}f\rangle . 
\end{equation*}
On the other hand, $[A, \mathcal{D}f] = \mathcal{D}
\langle A, \mathcal{D}f\rangle $. Indeed, we have for all $C \in C^\infty(E)$
\begin{equation*}
	\begin{aligned}
		\pi(A) \langle \mathcal{D}f, C\rangle  &= 2 \langle  \mathcal{D} \langle \mathcal{D}f, C\rangle , A\rangle , & \\
		\langle [A, \mathcal{D}f] + \mathcal{D} \langle A,
		\mathcal{D}f\rangle , C\rangle  + & \\ + \langle \mathcal{D}f,
		[A,C] + \mathcal{D}\langle A, C\rangle  \rangle  &= 2 \langle
		\mathcal{D}
		\langle \mathcal{D}f, C\rangle , A\rangle , &\quad  \text{by (5)}, \\ 
		\langle [A,\mathcal{D}f] + \mathcal{D}\langle A,
		\mathcal{D}f\rangle , C\rangle  + & \\ +  \langle \mathcal{D}f,
		[A,C]\rangle  &= 2 \langle \mathcal{D} \langle \mathcal{D}f, C\rangle , A\rangle ,
		& \quad \text{ by (4)}, \\
		\langle [A,\mathcal{D}f] + \mathcal{D}\langle A,\mathcal{D}f\rangle , C\rangle  + \frac{1}{2} \pi[A, C] f &=
		\frac{1}{2} \pi(A) \pi(C) f, & \\ \langle [A,\mathcal{D}f] +
		\mathcal{D}\langle A, \mathcal{D}f\rangle , C\rangle 
		&= \frac{1}{2} \pi(C) \pi(A) f, & \\ \langle [A, \mathcal{D}f] +
		\mathcal{D}\langle A, \mathcal{D}f\rangle ,
		C\rangle  &= 2 \langle \mathcal{D} \langle \mathcal{D}f, A\rangle , C\rangle,
	\end{aligned}
\end{equation*}
therefore we find
\begin{equation*}
	{[A}_\Lambda \mathcal{D} f] = 2 (\chi + \mathcal{D}) \langle A,
	\mathcal{D}f\rangle  = (\chi + \mathcal{D}) \pi(A) f.
\end{equation*}
We obtain then 
$	[A_\Lambda (S - \mathcal{D})f] = (S - \mathcal{D}) \pi(A)f,$
which vanish when imposing $S f = \mathcal{D}f, \forall f \in C^\infty(M)$. 

Note that we have also shown sesquilinearity of the form:
\begin{equation*}
	[A_\Lambda Sf] = (\chi + S) [A_\Lambda f].
\end{equation*}
It follows from (4) that $[Sf_\Lambda g] = [\mathcal{D}f_\Lambda g] = 0$, hence the
Lambda bracket (\ref{eq:lambdadef}) satisfies sesquilinearity.

Skew-symmetry is clear and we only need to check the Jacobi identity.
First let us compute $[A_\Lambda [B_\Gamma C]]$ for $A, B, C \in
C^\infty(\Pi E)$. By definition this is given by
\begin{multline*}
	[A_\Lambda [B, C] + (2 \eta + \mathcal{D}) \langle B, C\rangle ]  = [A, [B,C]] + (2 \chi
	+ \mathcal{D}) \langle A, [B,C]\rangle  + \\ + (2 \eta + \chi + \mathcal{D}) \pi(A) \langle
	B,C\rangle.
\end{multline*}
Similarly we find
\begin{equation*}
	[B_\Gamma [A_\Lambda C]] = [B, [A, C]] + (2 \eta + \mathcal{D}) \langle B, [A,C]\rangle 
	+ (2 \chi + \eta + \mathcal{D}) \pi(B) \langle A, C\rangle,
\end{equation*}
and finally
\begin{multline*}
	[ [A_\Lambda B]_{\Lambda + \Gamma} C] = [ [A,B]_{\Lambda +
	\Gamma} C] +  (\eta -  \chi) [\langle A, B\rangle _{\Lambda + \Gamma} C] = \\
	=[
	[A, B], C] + (2 \chi + 2 \eta + \mathcal{D}) \langle [A,B], C\rangle  + (\eta - \chi)
	\pi(C) \langle A, B\rangle.
\end{multline*}
Adding these last two equations and substracting the first,  we obtain the following. For the
coefficient of $\chi$ we find:
\begin{equation*}
	2\pi(B) \langle A, C\rangle  + 2 \langle [A,B], C\rangle  -  \pi(C) \langle A, B\rangle  - 2 \langle A, [B,C]\rangle  -
	\pi(A)\langle B,C\rangle,
\end{equation*}
and using the definition of the operator $\mathcal{D}$ this equals:
\begin{equation*}
	2 \Bigl( \pi(B) \langle A,C\rangle  + \langle [A,B], C\rangle  - \langle
	\mathcal{D} \langle A,B\rangle , C\rangle  
	- \langle A, [B,C]\rangle  - \langle \mathcal{D}\langle B,C\rangle , A\rangle
	\Bigr),
\end{equation*}
and this expression vanishes by (5). Similarly, the coefficient of
$\eta$ is given by:
\begin{equation*}
	2\langle B, [A, C]\rangle  + \pi(B) \langle A, C\rangle  + 2 \langle [A,B], C\rangle  + \pi(C) \langle A, B\rangle  - 2
	\pi(A) \langle B, C\rangle,
\end{equation*}
which in turn equals
\begin{equation*}
	2 \Bigl( \langle B, [A,C] + \mathcal{D}\langle A,C\rangle  \rangle  +
	\langle [A,B] + \mathcal{D}\langle A.B\rangle , C\rangle  -  \pi(A) \langle B,
	C\rangle  \Bigr),
\end{equation*}
and this vanishes by (5). Finally, the constant coefficient equals:
\begin{multline}
	Jac (A, B, C) + \mathcal{D} \langle  [A, B], C\rangle   + \mathcal{D}
	\langle B, [A, C]\rangle  + \mathcal{D} \pi(B) \langle A, C\rangle 
	\\ -  \mathcal{D} \langle A, [B, C]\rangle  - \mathcal{D} \pi(A) \langle B, C\rangle.
	\label{eq:constant1}
\end{multline}
\begin{lem}
	The following equality holds:
	\begin{equation*}
		\pi(B) \langle A, C\rangle  - \pi(A) \langle B, C\rangle  =
		\frac{4}{3} \langle [B, A], C\rangle +
	 \frac{2}{3} \langle A, [B, C]\rangle  + \frac{2}{3} \langle B,
	 [C,A]\rangle.
	\end{equation*}
	\label{lem:const.1}
\end{lem}
\begin{proof}
	This follows easily by applying (5) on the right hand side to
	obtain
	\begin{multline*}
		\pi(B) \langle A, C\rangle  - \pi(A) \langle B, C\rangle  =  2 \langle [B, A], C\rangle  \\ + \langle A,
		[B,C]\rangle  + \langle B, [C, A]\rangle  + \langle A, \mathcal{D}\langle
		B,C\rangle \rangle  - \langle B, \mathcal{D} \langle A,C\rangle \rangle.
	\end{multline*}
	Now using the definition of $\mathcal{D}$ in the last two terms the Lemma
	follows.
\end{proof}
replacing with the Lemma the corresponding terms in (\ref{eq:constant1})
we obtain that the constant term in the Jacobiator is:
\begin{equation*}
	\mathrm{Jac}(A,B,C) - \frac{1}{3} \mathcal{D} \langle [A,B], C\rangle  +
	\frac{1}{3} \mathcal{D} \langle B, [A,C]\rangle  - \frac{1}{3} \langle A,
	 [B,C]\rangle,
\end{equation*}
and by (2) in the definition of Courant Algebroid, we see that this vanishes.

We now need to check the Jacobi identity when one of the three terms is a
function $f \in C^\infty(M)$. So we compute for $A, B, \in
C^\infty(E)$:
\begin{equation*}
	[A_\Lambda [B_\Gamma f]] - [B_\Gamma [A_\Lambda f]] - [
	[A_\Lambda B]_{\Lambda + \Gamma} f] =  \pi(A) \pi(B) f - \pi(B)
	\pi(A) f - \pi [A, B] f,
\end{equation*}
and this vanishes by (1). 

All the other cases for the Jacobi identity are straightforward to check.
\end{proof}

\begin{rem}
It is clear that the construction above can be carried out locally and in a
way compatible with restriction maps, we have therefore constructed a sheaf of
$N_K=1$ SUSY Lie conformal algebras\footnote{Given the notation of \cite{bressler1}
and \cite{gerbes2}, one would be tempted to call these sheaves, SUSY vertex
algebroids. On the other hand, our bundles are of infinite rank, and they
correspond to the conformal Lie algebra associated to a vertex algebroid in the
usual case of op. cit.} associated to
any Courant algebroid $(E, \langle ,\rangle , [,], \pi) $. 
\end{rem}

The proof of the following is analogous to \cite[Thm 5.3]{heluani2}
\begin{prop}
	Let $(E, \langle,\rangle, [,], \pi)$ be a Courant algebroid, and
	$\cR$ be the corresponding sheaf of $N_K=1$ SUSY conformal Lie algebras
	constructed above. Let $\cU$ be the universal envelopping $N_K=1$ SUSY
	vertex algebra associated to $\cR$ \cite{heluani3}, we define
	$\Omega^\mathrm{ch}_M(E)$ to be the quotient of $\cU$ by the ideal
	generated by the relations ($f,g \in C^\infty(M)$, $A \in C^\infty (\Pi
	E)$, and $1_M$ is the constant function $1$)
	\begin{equation}
		:fg:=fg, \qquad :f A: = fA, \qquad 1_M = \vac. 
		\label{eq:ideal}
	\end{equation}
	Then $\Omega^\mathrm{ch}_M(E)$ is a sheaf of $N_K=1$ SUSY vertex algebras.
	When $E$ is the standard Courant algebroid of Example \ref{ex:courant1},
	$\Omega^\mathrm{ch}_M(E)$ is the chiral de Rham complex of $M$.
	\label{prop:universal}
\end{prop}

\begin{rem}
	Given a (SUSY)vertex algebra $V$ and a subset $I \subset V$, 
	in general it is difficult to compute the (SUSY) vertex algebra ideal
	generated by I. Indeed, one has to compute all products $v_{(n)} a$ for $v
	\in V$,  $a \in I$ and $n \in \mathbb{Z}$. This in particular includes all
	the OPEs (ie. $n \geq 0$) of the fields $v$ and $a$.

	The situation in Proposition \ref{prop:universal} is greatly simplified
	with the aid of (3) in the definition of Courant algebroids. Indeed,  
	 by the definition of the $\Lambda$ bracket on $\cR$ we have:
	\begin{equation}
		[A_\Lambda fB] = [A, fB] + (2 \chi + \mathcal{D}) f \langle A, B\rangle.
		\label{eq:wick1}
	\end{equation}
	On the other hand, the non-commutative Wick-formula implies: 
	\begin{equation}
		[A_\Lambda :fB:] = :(\pi(A)f) B: + :f[A,B]: + :f(2\chi + \mathcal{D})\langle
		A,B\rangle: 
		\label{eq:wick2}
	\end{equation}
	Now using (3) and the fact that $\mathcal{D}$ satisfies a Leibniz rule
	(e.g. $\mathcal{D}(fg)
	= f \mathcal{D} + \mathcal{D}(f)g$) we obtain substracting (\ref{eq:wick1}) and
	(\ref{eq:wick2}):
\begin{multline*}
	[A_\Lambda :fB:-fB] = \bigl(:(\pi(A)f) B: - (\pi(A)f)B\bigr) +
	\bigl(:f[A,B]: - f[A,B]\bigr) +
	\\ +
	2 \chi \bigl(:f \langle
	A,B\rangle: - f\langle A,B\rangle\bigr) + \bigl( :f \mathcal{D}\langle A,
	B\rangle : - f \mathcal{D}\langle A, B\rangle \bigr)
\end{multline*}
	which is a linear combination of the generators (\ref{eq:ideal}) of the
	ideal.
	\label{rem:1}
\end{rem}
\begin{rem}
	The explicit commutation relations of Proposition
	\ref{prop:susycommutators} are so simple that one is readily able to find
	relations between structures in a Courant algebroid $E$ and the
	corresponding sheaf of vertex algebras $\Omega^\mathrm{ch}_M(E)$. As an example, it is clear that
	given an integrable Dirac structure (c.f. \cite{gualtieri1} for a
	definition) we obtain a subalgebra of  $\Omega^\mathrm{ch}_M(E)$. Moreover,
	given a generalized complex structure $J$ on $T \oplus T^*$, we obtain 
	  two subsheaves of commutative vertex algebras inside the chiral de Rham
	  complex of $M$. This was
	showed for example in \cite{alekseev} and more recently in
	\cite{guttenberg} with relation to superfield brackets. 
	\label{rem:someremarkalpedo}
\end{rem}

\section{Another $N=2$ for K\"ahler manifolds}
\label{sec:anothern=2}
From now on, we will consider the standard Courant algebroid of Example
\ref{ex:courant1}, i.e. $E= T\oplus T^*$. We will use the same notation as in
\cite{heluani2}, in particular, for a coordinate system $\{x_i\}$ on $M$, we will
have the associated super-fields $\{B^i\}$ and $\{\Psi_i\}$ as in Example
\ref{ex:2.10}.

Given the above formalism, it is natural to expect extra symmetries of the chiral
de Rham complex, associated to extra structures on $T \oplus T^*$. Even though we
do not strictly need this notation, we will set up the formalism for a future
article relating the chiral de Rham complex with generalized complex geometries. Let $K \in
\fs\fo (T\oplus T^*)$. It can be written in the form:
\begin{equation*}
	K = \begin{pmatrix}
		I & \beta \\ B & -I^*	
	\end{pmatrix}
\end{equation*}
where $I={I_i}^j \in \End(T)$, $\beta = \beta^{ij} \in \Lambda^2 T$ is a bi-vector and
$B=B_{ij} \in \Lambda^2 T^*$ is a $2$-form.
The following is as easy generalization of 
\cite[Lem. 7.2]{heluani2}
\begin{lem}
	The assignment\footnote{From now on we omit the simbols $: :$ for normally
	ordered products when no confusion should arise.}
	\begin{equation*}
		K \mapsto J = ({I_i}^j SB^i) \Psi_j + \frac{1}{2} \bigl( \beta^{ij}
		\Psi_i \Psi_j + B_{ij} SB^i SB^j \bigr) + \Gamma^i_{jk} {I_i}^j
		TB^k,
	\end{equation*}
	defines a linear morphism
	\begin{equation*}
		\Gamma(M, \fs\fo(T \oplus T^*)) \rightarrow \Gamma(M,
		\Omega^\mathrm{ch}_M).
	\end{equation*}
	\label{lem:linearmorphism}
\end{lem}
Given a K\"ahler manifold $(M,g,J)$ with associated K\"ahler form $\omega$, we have two
commuting generalized complex structures (cf. \cite{gualtieri1}):
\begin{equation}
	J_1 =\begin{pmatrix}
		J & 0 \\ 0 & -J^*
	\end{pmatrix}, \qquad J_2 =\begin{pmatrix}
		0 & - \omega^{-1} \\ \omega & 0
	\end{pmatrix}.\label{eq:generalized}
\end{equation}
Recall \cite{heluani2} that the metric $g$ gives rise to a superconformal $N=1$
vector of central charge $3 \dim_\mathbb{R} M$:
\begin{equation}
	H = SB^i S\Psi_i + TB^i \Psi_i - TS \mathbf{g} = H_0 - TS \mathbf{g}.
	\label{eq:n=1def}
\end{equation}
where $\mathbf{g} = \log \sqrt{\det g_{ij}}$ and $H_0$ is the $N=1$ superconformal
vector of central charge $3 \dim_\mathbb{R} M$ constructed in \cite{malikov} (after twisting
to obtain a non-zero central charge). This $N=1$ structure is extended to
an $N=2$ structure by the superfield corresponding to $J_1$ if and only if $M$ is
Calabi-Yau \cite[Thm 7.4]{heluani2}. 

The symplectic case on the other hand is much less restrictive:
\begin{lem}
	Let $(M, \omega)$ be a symplectic $2n$-manifold with symplectic form $\omega =
	\omega_{ij}$. Let 
	\begin{equation}
		J = \frac{1}{2} \sum_{i,j = 1}^{2n} \left( (- \omega^{-1})_{ij} \Psi_i\Psi_j + \omega_{ij}SB^i
		SB^j \right),
		\label{eq:definoJ2}
	\end{equation}
	be the corresponding section of the chiral de Rham complex.
	Then the pair $\{H_0,J\}$ generates an $N=2$ vertex
	algebra of central charge $6n$.
	\label{lem:n=2}
\end{lem}
\begin{proof}
	We can compute the Lambda-brackets in any coordinate system, in particular,
	we can choose Darboux coordinates, to realize $\omega$ and $\omega^{-1}$ as
	constant matrices (locally). The Lemma is now reduced to a
	straightforward computation.
\end{proof}

On the other hand, if we are given a K\"ahler manifold $(M, g, \omega)$, the $N=2$
structure constructed by Lemma \ref{lem:n=2} does not extend the $N=1$ structure
associated to the Riemmanian manifold $(M,g)$. If we want to include the metric in
our $N=2$ structure, the 
 situation is much more subtle that in
the purely symplectic case. In particular we can now use the metric to raise or
lower indexes in $\omega$, but we cannot use Darboux coordinates unless $M$ is
flat. For a K\"ahler form $\omega=\omega_{ij}$, we will use $\omega^{ij} =
g^{ik}g^{jl} \omega_{kl}$. We have the following:
\begin{thm}
	Let $(M, J, g)$ be a Kahler $2n$-manifold with Kahler form $\omega$. Let
	$H$ be defined as in (\ref{eq:n=1def}). 
	Defining
	\begin{equation*}
		J = \frac{1}{2} \left( \omega_{ij} SB^i SB^j + \omega^{ij}\Psi_i
		\Psi_j \right), 
	\end{equation*}
	we obtain that the pair $\{J, H\}$ generates an $N=2$ vertex algebra of
	central charge $c = 6n$.
	\label{thm:anothern=2}
\end{thm}
\begin{proof}
	The proof can be found in the Appendix.
\end{proof}
\begin{rem}
	It follows from the proof of the Theorem that in fact, we have a
	\emph{family} of $N=2$ structures. Indeed, the superfield
	$H$, together with 
	\begin{equation*}
		J_\mu := \frac{1}{2} \left( \mu \, \omega_{ij} SB^i SB^j + \frac{1}{\mu}
		\omega^{ij} \Psi_i \Psi_j \right), \qquad \mu \in \mathbb{C}^\times,
	\end{equation*}
	generate an $N=2$ super vertex algebra of central charge $c=6n$.
	\label{rem:family}
\end{rem}
\section{Calabi-Yau case: $N=2,2$}\label{sec:n=2,2}

When $(M,J,g)$ is Calabi-Yau, we have two different $N=2$ structures sharing the
same underlying $N=1$ structure generated by (\ref{eq:n=1def}). We want to study
now how is that these two structures are related. For this we will need the
following Lemma, the proof of which can be found in the appendix:
\begin{lem}
	Let $J_i$, $i=1,2$ be the two superfields corresponding, by Lemma
	\ref{lem:linearmorphism}, to the two generalized complex structures
	(\ref{eq:generalized}) in the Calabi-Yau case and let $H$ be defined as
	in (\ref{eq:n=1def}), so that each pair $\{J_i, H\}$ generates an $N=2$
	vertex algebra of central charge $c = 6n$. 
	We have
	\begin{equation*}
		{[J_1}_\Lambda J_2] = {[J_1}_\Lambda J_2]|_{\Lambda = 0}.
	\end{equation*}
	\label{lem:secondcommutation}
\end{lem}
This Lemma allows one to define a section of CDR as $H':= - {[J_1}_\Lambda
J_2]$. For an explicit description of $H'$ in holomorphic coordinates see
(\ref{eq:htildedef}).

One would like to study if the algebra generated by $H$, $J_1$, $J_2$ and $H'$
closes, but the first term of (\ref{eq:htildedef}) makes the task of computing commutation
relations in this algebra a tedious task. We can avoid lots of explicit
computations by using the Jacobi identity for SUSY $N_K=1$ Lie conformal algebras.
This is illustrated in the proof of the following theorem, which can be found in the
Appendix. 
\begin{thm}
	Let $(M,g,J)$ be a Calabi-Yau $2n$-manifold with K\"ahler form $\omega$. Let
	$H,\, J_1,\, J_2$ and $H'$ be the global sections of the chiral de Rham
	complex of $M$ defined above. 
	Define 
	\begin{equation}
		H^\pm = \frac{1}{2} (H \pm H'), \qquad J^\pm =
		\frac{1}{2} (J_1 \pm J_2).
		\label{eq:n=2,2def}
	\end{equation}
	Then each pair $(J^+, H^+)$ and $(J^-, H^-)$ generates an $N=2$ vertex algebra of
	central charge $c =
	3n$. Moreover, these two different $N=2$ structures commute, namely, the
	superfields $\{J^\pm, H^\pm\}$ generate the tensor product of two $N=2$
	vertex algebras of central charge $c=3n$.
	\label{thm:n=2,2}
\end{thm}
\begin{rem}
	Note that neither of the $N=2$ structures is purely holomorphic or
	anti-holomorphic. On the other hand, the central charge of each sector
	$\{H^\pm, J^\pm\}$ agrees with the central charge of the holomorphic chiral
	de Rham complex of $M$ (after untwisting). 
	\label{rem:central charges}
\end{rem}
Even though the metric was not involved in defining $\Omega^\mathrm{ch}_M(E)$, the
existence of these extra symmetries in the Calabi-Yau case is intimately dependent
on the metric. Moreover, recall from \cite{heluani2} that the existence of an $N=1$
supersymmetry was enough for us to define superfields in terms of the original
fields of \cite{malikov}. Recall (Example \ref{ex:2.10}) that these expressions
are of the form:
\begin{equation*}
	B^i(Z) = b^i(z) + \theta \phi^i(z), \qquad \Psi_i(Z) = \psi_i(z) + \theta
	a_i(z),
\end{equation*}
where $\{b^i, a_i, \psi_i, \phi^i\}$ are generators of the usual $bc-\beta\gamma$
system.

We can start however with $H'$ as the generator of our supersymmetry as follows.
Define the operator $S' = H'_{(0|1)}$, that is the coefficient of $z^{-1}$ in the
superfield $H'(z,\theta)$. One can compute the commutation relation of $H'$ with
itself as in (\ref{eq:nosecual5}), from where it follows that $(S')^2 = T$. Note
however that $S' \neq S$. In particular $S'$ is not the zero mode of an $N=1$
superconformal vector (c.f. \cite{heluani3} for a definition). Recall from
\cite{malikov} that the
fields $\phi^i$ (resp. $\psi_i$) transform as $1$-forms on $M$ do (resp. vector
fields), we can view the
metric $g_{ij}$ as an isomorphism $T\simeq T^*$, and as before use it to lower and
raise indexes. The interplay of this isomorphism and the supersymmetry generated by
$S'$ is explained in the
following
\begin{thm}
	Let $(M,g,J)$ be a Calabi-Yau manifold.
	The assignement
	\begin{equation}
		\begin{aligned}
			\psi_i &\mapsto \phi_i:=g_{ij} \phi^j, &\qquad \phi^i
			&\mapsto g^{ij}\psi_{j}=S'b^i, \\
			b^i &\mapsto b^i, &\qquad a_i &\mapsto S'\phi_i
		\end{aligned}
		\label{eq:isomorphismstrange}
	\end{equation}
	defines an automorphism of the chiral de Rham complex of $M$ as a sheaf of
	vertex algebras.
	\label{thm:isomorphismstrange}
\end{thm}
\begin{rem}
	Note that we have computed explicitly $S'b^i$, while we
	didn't compute $S' \phi_i$. This latter field is rather complicated and we
	don't need its specific description to state this Theorem (see however
	(\ref{eq:sprimeonphi}) below). On the other
	hand, we could have phrased this as a statement of SUSY vertex algebras
	rather than vertex algebras, namely,
	 if one defines the superfields
	\begin{equation*}
		A^i(z,\theta) = b^i(z) + \theta S'b^i(z), \qquad
		\Theta_i(z,\theta) = \phi_i(z) + \theta S'\phi_i(z)
	\end{equation*}
	Then the assignement $B^i \mapsto A^i$ and $\Psi_i \mapsto \Theta_i$
	preserves Lambda brackets (ie. OPEs). But note that it doesn't preserve the
	differential $S$ (it is rather mapped to $S'$), hence it is not an automorphism of SUSY vertex algebras.
	Indeed, these are two different  $N_K=1$ SUSY vertex algebra
	structures on $\Omega_M^\mathrm{ch}$. 
	\label{rem:isomorphismstrange}
\end{rem}
\begin{rem}
	Note that in the flat case, this automorphism reduces to the identity
	automorphism on the bosonic part:
	\begin{equation*}
		\begin{aligned}
			b^i &\mapsto b^i, &\qquad a_i &\mapsto a_i, \\
			\psi_i &\mapsto g_{ij}\phi^j, &\qquad \phi^i &\mapsto
			g^{ij} \psi_j,
		\end{aligned}
	\end{equation*}
	where $g_{ij}$ and $g^{ij}$ are inverse constant matrices.
	\label{rem:isomosphism2}
\end{rem}
\begin{proof}[Proof of Theorem \ref{thm:isomorphismstrange}] \hfill 
	
	By a straightforward
	computation we find using (\ref{eq:htildedef}) 
	\begin{equation}
		\begin{aligned}
			S'\phi^k &= \left( \Gamma^{k}_{jl} g^{ij} \phi^l \right) \psi_i -
		\left( \Gamma^m_{jl} g^{kj} \phi^l \right)\psi_m + g^{ik}a_i \\
		S' g_{jk} &= g^{il} g_{jk,i}\psi_l.
	\end{aligned}
		\label{eq:sprimeonphi}
	\end{equation}
	From where it follows, using that the metric is parallel with respect to
	the Levi-Civita connection:
	\begin{equation}
		{[\phi^i}_\lambda S'\phi_j] = \Gamma^i_{jl} \phi^l + \Gamma^k_{mj}
		g^{im} \phi_k.
		\label{eq:nosecualesahora2}
	\end{equation}
	Note that we use lambda brackets instead of SUSY Lambda brackets because we
	are only interested in a vertex algebra isomorphism.
	Similarly, we find
	${[g_{pi}}_\lambda S' \phi_j] = - g_{pi,j}$,
	and combining with (\ref{eq:nosecualesahora2}) this gives:
	\begin{multline*}
		{[S'\phi_j}_\lambda g_{pi} \phi^i] = g_{pi,j} \phi^i -
		\Gamma^{i}_{jl} g_{pi} \phi^l - \Gamma^{k}_{mj} g_{pi} g^{im}
		\phi_k = \\ = \Gamma^l_{jp} \phi_l + \Gamma^l_{ji} g_{pl} \phi^i -
		\Gamma^i_{jl} g_{pi} \phi^l - \Gamma^k_{pj} \phi_k
		\label{eq:nosecualesahora4}
	\end{multline*}
	And this implies clearly
	\begin{equation}
		{[S'\phi_j}_\lambda \phi_i] = {[\phi_i}_\lambda S'\phi_j] = 0.
		\label{eq:nosecualesahora5}
	\end{equation}
	Using that $S'$ squares to $T$ and that the zero mode of a field is a
	derivation of all $n$-th products on a vertex algebra, this last equation
	gives:
	\begin{equation}
		{[S'\phi_i}_\lambda S'\phi_j] = 0.
		\label{eq:nosecualesahora6}
	\end{equation}
	From (\ref{eq:sprimeonphi}) we easily find:
	\begin{multline}
		{[b^i}_\lambda S'\phi_j] = {[b^i}_\lambda (S' g_{jk}) \phi^k] +
		{[b^i}_\lambda g_{jk} S'\phi^k] = \\ = 0 + g_{jk} {[b^i}_\lambda
		g^{lk}a_l] = - g_{jk} g^{ik} = -{\delta^{i}}_j.
		\label{eq:nosecualesahora7}
	\end{multline}
	and again using that $S'$ is a derivation of the bracket that squares to
	$T$, we obtain 
	\begin{equation}
		{[S'b^i}_\lambda S'\phi_j] = 0.
		\label{eq:nosecualesahora8}
	\end{equation}
	Equations (\ref{eq:nosecualesahora5})-(\ref{eq:nosecualesahora8}) easily
	imply that the only non-vanishing lambda brackets among the fields
	$\{b^i,\, S'b^i,\, \phi_i,\, S'\phi_i\}$ are:
	\begin{equation*}
		{[b^i}_\lambda S'\phi_j] = - {\delta^i}_j, \qquad {[\phi_i}_\lambda
		S'b^j] = {\delta_i}^j,
	\end{equation*}
	proving the Theorem.
\end{proof}
Let us study how this automorphism acts on the structures described so far in a
Calabi-Yau manifold. Since we have an automorphism of vertex algebras which does
not preserve the SUSY structure, we need to work with usual fields (as opposed to
superfields). We can describe the $N=2,2$ structure defined above in terms of the
usual fields $\{a_i, b^i, \psi_i, \phi^i\}$ of
\cite{malikov} plus the supersymmetry generator $S$.
We define\footnote{Here the expression for $H'$ is valid only in holomorphic
coordinates. In general, it can be defined as the $0$-th product of the fields $SJ_1$ and
$J_2$.}:
\begin{equation*}
	\begin{aligned}
	J_1 &= \left({\omega_i}^j \phi^i \right) \psi^j +
	\Gamma^i_{jk} {\omega_i}^j Tb^k , \\
	J_2 &=  \frac{1}{2} \left( \omega^{ij}\psi_i\psi_j +
	\omega_{ij} \phi^i \phi^j \right), \\
	H &= \phi^i a_i + Tb^i \psi_i - T(\mathbf{g}_{,i} \phi^i), \\
	H'&= \Gamma^{k}_{jl} g^{ij} \phi^l (\psi_i \psi_k) + g^{ij}a_i \psi_j +
	g_{ij} Tb^i \phi^j,\\
	H^\pm &= \frac{1}{2} (H \pm H'), \qquad J^\pm = \frac{1}{2} (J_1 \pm J_2), \qquad S = H_{(0)}.
\end{aligned}
\end{equation*}
Then Theorem \ref{thm:n=2,2} says that for a Calabi Yau $2n$-manifold, the fields
\[\{H^\pm, J^\pm, SJ^\pm,
SH^\pm\},\] generate two commuting copies of the $N=2$ super vertex algebra of
central charge $c = 3n$.

It is easy to show that the fields $J_i$ are invariant under the automorphism of
Theorem \ref{thm:isomorphismstrange}. And it follows from (\ref{eq:nosecual3}) and
(\ref{eq:nosecual4}) that this automorphism exchanges $SJ_1 \leftrightarrow SJ_2$.
Similarly it follows from (\ref{eq:htildedef}) that this automorphism exchanges $H
\leftrightarrow H'$. Finally, it follows from (\ref{eq:nosecual5}) and
(\ref{eq:nosecual6}) that $SH$ and $SH'$ are invariant by this automorphism.
Therefore we have:
\begin{prop}
	For a Calabi-Yau manifold $(M,g,J)$, the $N=2$ subalgebra of CDR generated
	by $\{J^+, H^+, SJ^+, SH^+\}$ is invariant under the automorphism of Theorem
	\ref{thm:isomorphismstrange}. On the other hand, this automorphism, when
	restricted to the $N=2$ subalgebra generated by $\{J^-, H^-, SJ^-, SH^-\}$
	is the automorphism which is the identity in the even part and
	multiplication by $-1$ in the odd part. 

	Moreover, the original $N=2$ subalgebra of \cite{heluani2} generated by
	$\{J_1, H, SJ_1, SH\}$ is mapped to the vertex algebra generated by $\{J_1,
	H', SJ_2, SH'\}$ therefore this algebra is another $N=2$ of central charge
	$c = 6n$. Similarly, the $N=2$ vertex algebra of Theorem
	\ref{thm:anothern=2} generated by $\{J_2, H, SJ_2, SH\}$ is mapped to the
	vertex algebra generated by $\{J_2, H', SJ_1, SH\}$ and therefore this
	algebra is another $N=2$ vertex algebra of central charge $c = 6n$.
	\label{prop:automorphimaction}
\end{prop}
\begin{rem}
	Note that if we use $\mu = \sqrt{-1}$ as in Remark \ref{rem:family} in
	the definition of $J_2$, then the automorphism of Theorem
	\ref{thm:isomorphismstrange} maps the corresponding $U(1)$ current $J_2$ to
	$-J_2$, which looks more like the \emph{mirror involution} of the $N=2$
	super vertex algebra. We plan to return to these matters in the future.
	\label{rem:future1}
\end{rem}
\begin{rem}
	We note that the expression of $H'$ in (\ref{eq:htildedef}) does not seem
	to depend in the complex structure. However, we have used holomorphic
	coordinates to arrive to such expression. 
\end{rem}
\section{Hyper-K\"ahler case: $N=4,4$}\label{sec:n=4,4}
Let $(M, g, I, J, K)$ be a Hyper-K\"ahler manifold, that is a Riemannian manifold
$(M,g)$, together with three complex structures $I,\, J,\, K$ satisfying the
quaternionic relations
\begin{equation*}
	IJ = -JI = K,
\end{equation*}
and such that $(M, g, I)$, $(M,g,J)$ and $(M,g,K)$ are K\"ahler. Lemma
\ref{lem:linearmorphism} associates three superfields $J_i$, $i=1,2,3$ to these
complex structures and three superfields $J_{\omega_i}$ to the corresponding
K\"ahler forms. Define 
\begin{equation*}
	J_i^\pm = \frac{1}{2} \left(J_i \pm J_{\omega_i} \right),
\end{equation*}

It was shown in \cite[Thm. 7.4]{heluani2} that the superfields
$\{H, J_i\}$ generate an $N=4$ super vertex algebra of central charge $c = 3
\dim_\mathbb{R}
M$. We can extend this Theorem as follows.

\begin{thm}
Let $(M, g, I,J,K)$ be a Hyper-K\"ahler $4n$-manifold. 
	\begin{enumerate}
		\item The superfield $H'$ does not depend on the complex structure
			used, namely 
			\begin{equation*}
				{[J_1}_\Lambda J_{\omega_1}] =  {[J_2}_\Lambda
				J_{\omega_2}] = {[J_3}_\Lambda J_{\omega_3}].
			\end{equation*}
			Define then $H^\pm = \tfrac{1}{2} (H \pm H')$.
		\item
	The superfields
	$H$, $J_1$, $J_{\omega_2}$ and  $J_{\omega_3}$ generate an $N=4$ super-vertex algebra
	of central charge $c = 12n$.
		\item The fields $\{H^+, J_1^+, J_2^+, J_3^+\}$ and $\{H^-, J_1^-,
			J_2^-, J_3^-\}$ generate
			two commuting copies of the $N=4$ super vertex algebra of
			central charge $c = 6n$.
	\end{enumerate}
	\label{thm:n=4,4}
\end{thm}
\begin{proof}
	The proof of this theorem can be found in the Appendix.
\end{proof}

\appendix
\section{Appendix}
\begin{proof}[Proof of Theorem \ref{thm:anothern=2}] \hfill 

	We can compute the Lambda-brackets in any coordinate system, hence let us
	choose holomorphic coordinates. Define
	\begin{equation*}
		\beta =  \omega^{\alpha, \bar\beta} \Psi_\alpha
		\Psi_{\bar{\beta}}, \qquad \Omega = \omega_{\alpha, \bar{\beta}}
		SB^\alpha SB^{\bar{\beta}},
	\end{equation*}
	then we have $J = \beta + \Omega$. We first start by computing
	$[\beta_\Lambda \beta]$:
	\begin{equation*}
			{[\omega^{\gamma, \bar{\delta}}}_\Lambda \beta] =
			\omega^{\alpha, \bar{\beta}} {\omega^{\gamma,
			\bar{\delta}}}_{, \alpha} \Psi_{\bar{\beta}}
			- \omega^{\alpha,\bar{\beta}} {\omega^{\gamma,
			\bar{\delta}}}_{,\bar{\beta}} \Psi_{\alpha} ,  
	\end{equation*}
	Since $\omega$ is parallel, we have
	\begin{equation*}
		{[\omega^{\gamma,\bar{\delta}}}_\Lambda \beta] = -
		\Gamma^{\gamma}_{\varepsilon, \alpha} \omega^{\alpha,
		\bar{\beta}}\omega^{\varepsilon, \bar{\delta}} \Psi_{\bar{\beta}} 
		+ \Gamma^{\bar{\delta}}_{\bar{\varepsilon}\bar{\beta}}
		\omega^{\alpha\bar{\beta}}\omega^{\gamma\bar{\delta}} \Psi_\alpha = [\beta_\Lambda \omega^{\gamma,\bar{\delta}}].
	\end{equation*}
	Similarly:
	\begin{equation*}
		\begin{aligned}
			{[\Psi_\gamma}_\Lambda \beta] &= {\omega^{\alpha,
				\bar{\beta}}}_{,\gamma} \Psi_\alpha
				\Psi_{\bar{\beta}}  \\
				&= - \Gamma^\alpha_{\varepsilon,\gamma}
				\omega^{\varepsilon, \bar{\beta}} \Psi_\alpha
				\Psi_{\bar{\beta}} =
					[\beta_\Lambda \Psi_\gamma],\end{aligned}
	\end{equation*}
	hence
	\begin{equation*}
		[\beta_\Lambda \Psi_\gamma \Psi_{\bar{\delta}}] = -
		\Gamma^\alpha_{\varepsilon, \gamma} \omega^{\varepsilon,
		\bar{\beta}} \Psi_\alpha \Psi_{\bar{\beta}} \Psi_{\bar{\delta}} 
		+ \Gamma^{\bar{\beta}}_{\bar{\varepsilon}\bar{\delta}}
		\omega^{\alpha\bar{\varepsilon}} \Psi_{\gamma} \Psi_\alpha
		\Psi_{\bar{\beta}}.
	\end{equation*}
	Finally we get
	\begin{equation*}
		\begin{aligned}
			{[\beta}_\Lambda \beta] &= -
		\Gamma^\gamma_{\varepsilon,\alpha} \omega^{\alpha, \bar{\beta}}
		\omega^{\varepsilon, \bar{\delta}}
		\Psi_{\bar{\beta}} \Psi_\gamma \Psi_{\bar{\delta}}  -
		\Gamma^\alpha_{\varepsilon, \gamma} \omega^{\gamma, \bar{\delta}}
		\omega^{\varepsilon, \bar{\beta}} \Psi_\alpha
		\Psi_{\bar{\beta}}\Psi_{\bar{\delta}} + c.c. \\
		&= - 2 \Gamma^{\gamma}_{\varepsilon, \alpha} \omega^{\alpha,
		\bar{\beta}} \omega^{\varepsilon, \bar{\delta}}
		\Psi_{\bar{\beta}}\Psi_{\gamma} \Psi_{\bar{\delta}} + c.c.
		\end{aligned}
	\end{equation*}
	where $c.c.$ denotes the complex conjugate.
	This in turn can be expressed as:
	\begin{equation*}
		\begin{aligned}
			{[\beta}_\Lambda\beta] &= 2 \Gamma^{\gamma}_{\varepsilon,
			\alpha} \omega^{\alpha, \bar{\beta}} \omega^{\varepsilon,
			\bar{\delta}} \Psi_{\bar{\delta}} \Psi_{\gamma}
			\Psi_{\bar{\beta}} + c.c \\ &= 2 \Gamma^{\gamma}_{\alpha,
			\varepsilon} \omega^{\alpha, \bar{\delta}}
			\omega^{\varepsilon, \bar{\beta}}
			\Psi_{\bar{\beta}}\Psi_{\gamma}\Psi_{\bar{\delta}} +
			c.c.\\&= 2\Gamma^{\gamma}_{\varepsilon, \alpha} \omega^{\varepsilon,
			\bar{\delta}} \omega^{\alpha, \bar{\beta}}
			\Psi_{\bar{\beta}}\Psi_\gamma\Psi_{\bar{\delta}} + c.c. =
			-[\beta_\Lambda\beta].
		\end{aligned}
	\end{equation*}
	Hence we have $[\beta_\Lambda\beta] = [\Omega_\Lambda\Omega] = 0$ and we
	are left to compute $[\beta_\Lambda \Omega]$. For that we will need
	\begin{equation*}
		{[\Psi_\alpha}_\Lambda SB^\gamma SB^{\bar{\delta}}] = \chi
		\delta_{\alpha}^\gamma SB^{\bar{\delta}},
	\end{equation*}
	hence
	\begin{equation*}
		\begin{aligned}
			{[\Psi_\alpha}_\Lambda \Omega] &= 
		(\omega_{\gamma,\bar{\delta}})_{,\alpha} SB^{\gamma}SB^{\bar{\delta}} +
		 \chi \omega_{\alpha, \bar{\delta}} SB^{\bar{\delta}}\\ 
		&=  \Gamma^{\varepsilon}_{\gamma,\alpha}
		\omega_{\varepsilon, \bar{\delta}} SB^\gamma SB^{\bar{\delta}} +
		 \chi \omega_{\alpha, \bar{\delta}} SB^{\bar{\delta}}, 
	\end{aligned}
	\end{equation*}
	and by skew-symmetry we have:
	\begin{equation*}
		\begin{aligned}
			{[\Omega}_\Lambda \Psi_{\alpha}] &= 
			\Gamma^{\varepsilon}_{\gamma,\alpha}
			\omega_{\varepsilon,\bar{\delta}} SB^\gamma
			SB^{\bar{\delta}} -  \Gamma^\varepsilon_{\alpha,
			\gamma} \omega_{\varepsilon \bar{\delta}} SB^\gamma
			SB^{\bar{\delta}} - \\ & \quad -  
			\Gamma^{\bar{\varepsilon}}_{\bar{\delta}\bar{\gamma}}
			\omega_{\alpha,\bar{\varepsilon}} SB^{\bar{\gamma}}
			SB^{\bar{\delta}} - \omega_{\alpha
			\bar{\delta}} TB^{\bar{\delta}}  -  \chi
			\omega_{\alpha\bar{\delta}} SB^{\bar{\delta}} \\ &= -  
			\omega_{\alpha\bar{\delta}} TB^{\bar{\delta}}-  \chi
			\omega_{\alpha\bar{\delta}} SB^{\bar{\delta}}.
		\end{aligned}
	\end{equation*}
	And now we can compute then:
	\begin{equation*}
		\begin{aligned}
			{[\Omega}_\Lambda \Psi_\alpha \Psi_{\bar{\beta}}] &=
			- \left( \omega_{\alpha\bar{\delta}} TB^{\bar{\delta}}
			 + \chi
			\omega_{\alpha \bar{\delta}} SB^{\bar{\delta}} 
			\right) \Psi_{\bar{\beta}} -  \Psi_\alpha \left(
			\omega_{\delta \bar{\beta}} TB^\delta + \chi \omega_{\delta
			\bar{\beta}} SB^\delta\right) - \\ & \quad -
			 \int_0^\Lambda [{\omega_{\alpha \bar{\delta}}
			TB^{\bar{\delta}} + \chi \omega_{\alpha \bar{\delta}}
			SB^{\bar{\delta}}}_\Gamma \Psi_{\bar{\beta}}] d\Gamma.
		\end{aligned}
	\end{equation*}
	We use quasi-commutativity now to obtain:
	\begin{equation*}
		\begin{aligned}
			{[\Omega}_\Lambda \Psi_\alpha \Psi_{\bar{\beta}}]  &=
			- \left( \omega_{\alpha\bar{\delta}} TB^{\bar{\delta}}
			 + \chi
			\omega_{\alpha \bar{\delta}} SB^{\bar{\delta}} 
			\right) \Psi_{\bar{\beta}} + \frac{1}{2} \chi
			T(\omega_{\alpha, \bar{\beta}})  -
			\frac{1}{2} \chi \lambda \omega_{\alpha \bar{\beta}} - c.c.
		\end{aligned}
	\end{equation*}
	And this easily implies:
	\begin{equation}
		\begin{aligned}
			{[\Omega}_\Lambda \beta] &= - \omega^{\alpha
			\bar{\beta}} \Biggl[ \left( \omega_{\alpha \bar{\delta}}
			TB^{\bar{\delta}} + \chi \omega_{\alpha \bar{\delta}}
			SB^{\bar{\delta}} \right)  \Psi_{\bar{\beta}} \Biggr]+ \\ &
			\quad + \frac{1}{2} \chi \omega^{\alpha \bar{\beta}}
			T(\omega_{\alpha \bar{\beta}}) - \frac{1}{2} \chi \lambda
			\omega^{\alpha \bar{\beta}} \omega_{\alpha\bar{\beta}} +
			c.c. \\ &= - 
			TB^{\bar{\delta}}\Psi_{\bar{\delta}} -
			 \chi SB^{\bar{\delta}} \Psi_{\bar{\delta}} +
			 \chi T(\omega^{\alpha \bar{\beta}})
			\omega_{\alpha \bar{\beta}} + \\ & \quad + \frac{1}{2} \chi
			\omega^{\alpha \bar{\beta}}T(\omega_{\alpha \bar{\beta}}) -
			\frac{n}{2}  \lambda \chi +
			c. c. \\ &= - 
			TB^{\bar{\delta}}\Psi_{\bar{\delta}} -
			 \chi SB^{\bar{\delta}} \Psi_{\bar{\delta}} -
			\frac{1}{2} \chi T(\omega^{\alpha \bar{\beta}})
			\omega_{\alpha \bar{\beta}}  -
			\frac{n}{2} \lambda \chi +
			c. c.
		\end{aligned}
		\label{eq:proof2}
	\end{equation}
	The third term in this last expression is given by:
	\begin{multline*}
			T(\omega^{\alpha \bar{\beta}}) \omega_{\alpha
			\bar{\beta}} = {\omega^{\alpha \bar{\beta}}}_{,\gamma}
			\omega_{\alpha\bar{\beta}} TB^\gamma + c.c. = 
			\Gamma^{\alpha}_{\gamma \varepsilon} \omega^{\varepsilon
			\bar{\beta}}\omega_{\alpha\bar{\beta}} TB^\gamma + c.c. =
			\\ =
			\Gamma^\alpha_{\gamma \alpha} TB^\gamma + c.c. =
			\mathbf{g}_{,\gamma}
			TB^\gamma + c.c. = T \mathbf{g},
	\end{multline*}
	hence (\ref{eq:proof2}) reads:
	\begin{equation*}
		{[\Omega}_\Lambda \beta] = - TB^{i}\Psi_{i} - \chi
		SB^i \Psi_i - \chi T\mathbf{g} - n\lambda \chi.
	\end{equation*}
	And now using skew-symmetry we easily find:
	\begin{equation*}
		{[J}_\Lambda J] = - TB^i\Psi_i - SB^iS\Psi_i + TS\mathbf{g} - 2n
		\lambda \chi = - \left( H + \frac{c}{3} \lambda\chi \right)
	\end{equation*}
	where $c = 6n$. 

	We are left to check that $J$ is a primary field of conformal weight $1$
	\cite{heluani3}.
	
	For this we recall from \cite[(7.4.8)]{heluani2} that we have:
	\begin{equation*}
		\begin{aligned}
			{[H}_\Lambda SB^\alpha] &= (2T + \lambda + \chi S)
			SB^\alpha, \\ [H_\Lambda \Psi_\alpha] &= (2T + \lambda +
			\chi S) \Psi_\alpha + \lambda \chi \mathbf{g}_{,\alpha},\\
			[H_\Lambda \omega^{\alpha \bar{\beta}}] &= (2T + \chi S)
			\omega^{\alpha \bar{\beta}}.
		\end{aligned}
	\end{equation*}
	From this it follows easily that $\Omega$ is primary of conformal weight
	$1$. On the other hand, we obtain:
	\begin{equation*}
		{[H}_\Lambda \Psi_\alpha\Psi_{\bar{\beta}}] = (2T + 2\lambda +
		\chi S) \Psi_\alpha \Psi_{\bar{\beta}} + \lambda\chi
		\left(\mathbf{g}_{,\alpha} \Psi_{\bar{\beta}} -
		\mathbf{g}_{\bar{\beta}} \Psi_{\alpha} \right),
	\end{equation*}
	Therefore:
	\begin{multline}
		{[H}_\Lambda \beta] =  \Bigl( (2T + \chi S)
		\omega^{\alpha\bar{\beta}} \Bigr) \Psi_\alpha\Psi_{\bar{\beta}} +
		 \omega^{\alpha\bar{\beta}} (2T+2\lambda + \chi S) \Psi_{\alpha}
		\Psi_{\bar{\beta}} + \\ +  \lambda \chi \omega^{\alpha \bar{\beta}}
		(\mathbf{g}_{,\alpha} \Psi_{\bar{\beta}} -
		\mathbf{g}_{,\bar{\beta}} \Psi_\alpha) + 
		\int_0^\Lambda (-2\lambda -\chi \eta)
		{[\omega^{\alpha\bar{\beta}}}_\Lambda \Psi_\alpha
		\Psi_{\bar{\beta}}] d\Gamma =\\  =  (2T + 2\lambda +
		\chi S) \beta + \lambda \chi \omega^{\alpha \bar{\beta}}
		(\mathbf{g}_{, \alpha} \Psi_{\bar{\beta}} -
		\mathbf{g}_{,\bar{\beta}} \Psi_\alpha) + \\ +  \lambda\chi
		{\omega^{\alpha \bar{\beta}}}_{,\alpha} \Psi_{\bar{\beta}}
		- \lambda \chi {\omega^{\alpha
		\bar{\beta}}}_{,\bar{\beta}} \Psi_{\alpha},
		\label{eq:nosecuales}
	\end{multline}
	and using the fact that $\omega$ is parallel and that in a K\"ahler manifold
	$\Gamma^{\alpha}_{\alpha\gamma} = \mathbf{g}_{,\gamma}$, we see that this
	last expression equals:
	\begin{equation*}
		{[H}_\Lambda \beta] = (2T + 2\lambda +\chi S) \beta,
	\end{equation*}
	thus proving that $\beta$ and therefore $J$ is primary of conformal weight
	$1$.
\end{proof}
\begin{proof}[Proof of Lemma \ref{lem:secondcommutation}]\hfill 
	We will use holomorphic coordinates again and use the same notation as in
	the Proof of Theorem \ref{thm:anothern=2}. We need:
	\begin{equation}
			{[J_1}_\Lambda \Psi_\alpha] = -i (\chi + S) \Psi_\alpha - i
			\lambda \mathbf{g}_{,\alpha}. 
			\label{eq:prueba2}
	\end{equation}
	Hence we have
	\begin{equation*}
		\begin{aligned}
		{[J_1}_\Lambda \Psi_{\alpha}\Psi_{\bar{\beta}}] &= - i \left(
		(\chi+S)\Psi_\alpha
		\right) \Psi_{\bar{\beta}} - i \lambda \mathbf{g}_{,\alpha}
		\Psi_{\bar{\beta}} - i \Psi_\alpha (\chi + S) \Psi_{\bar{\beta}} -
		i \Psi_{\alpha} \lambda \mathbf{g}_{,\bar{\beta}}\\ &= - i
		(S\Psi_\alpha)\Psi_{\bar{\beta}} - i \Psi_{\alpha}
		S\Psi_{\bar{\beta}} - i \lambda \left( \mathbf{g}_{,\alpha}
		\Psi_{\bar{\beta}} + \mathbf{g}_{,\bar{\beta}} \Psi_{\alpha}
		\right).
	\end{aligned}
	\end{equation*}
	We also have:
	\begin{equation*}
		\begin{aligned}
		{[J_1}_\Lambda \omega^{\alpha \bar{\beta}}] &= - i {\omega^{\alpha
		\bar{\beta}}}_{,\gamma} SB^\gamma + i {\omega^{\alpha,
		\bar{\beta}}}_{,\bar{\gamma}} SB^{\bar{\gamma}}, \\ 
		&=  i \Gamma^{\alpha}_{\varepsilon\gamma} \omega^{\varepsilon
		\bar{\beta}} SB^\gamma - i
		\Gamma^{\bar{\beta}}_{\bar{\varepsilon}\bar{\gamma}}
		\omega^{\alpha\bar{\varepsilon}} SB^{\bar{\gamma}},
	\end{aligned}
	\end{equation*}
	from where it easily follows now:
	\begin{multline}
			{[J_1}_\Lambda \beta] = \left( i
			\Gamma^{\alpha}_{\varepsilon\gamma} \omega^{\varepsilon
			\bar{\beta}} SB^\gamma
			\right)(\Psi_{\alpha}\Psi_{\bar{\beta}}) -
			i\omega^{\alpha\bar{\beta}} \left(
			(S\Psi_\alpha)\Psi_{\bar{\beta}} \right) - \\ -  i\lambda \omega^{\alpha
			\bar{\beta}} (\mathbf{g}_{,\alpha} \Psi_{\bar{\beta}}) 
			+ i \int_0^\Lambda [\Gamma^\alpha_{\varepsilon\gamma}
			\omega^{\varepsilon \bar{\beta}} {SB^\gamma}_\Gamma
			\Psi_\alpha \Psi_{\bar{\beta}}] d\Gamma
			+ c.c.
			\label{eq:primeracomplicada}
	\end{multline}
	We recognize on the first two terms (plus their complex conjugates) the
	expression\footnote{note that there are no quasi-associativity issues in
	the second term of this expression, that is why we can write it without
	parenthesis.}:
	\begin{equation*}
		-\Gamma^{k}_{jl}g^{ij}SB^l (\Psi_i \Psi_k) - g^{ij} S\Psi_i
		\Psi_j.
	\end{equation*}
 	We therefore need to compute the integral term in
	(\ref{eq:primeracomplicada}):
	\begin{equation*}
		i \int_0^\Lambda \eta \left( \Gamma^{\alpha}_{\varepsilon\alpha}
		\omega^{\varepsilon \bar{\beta}} \right) \Psi_{\bar{\beta}} + c.c.
		=  i \lambda \omega^{\alpha \bar{\beta}} \mathbf{g}_{\alpha}
		\Psi_{\bar{\beta}} + c.c.
	\end{equation*}
	cancelling the third term in (\ref{eq:primeracomplicada}). We obtain then
	\begin{equation}
		{[J_1}_\Lambda \beta] = -\Gamma^{k}_{jl}g^{ij}SB^l (\Psi_i \Psi_k) - g^{ij} S\Psi_i
		\Psi_j.
		\label{eq:conbeta}	
	\end{equation}
	To compute ${[J_1}_\Lambda \Omega]$ we need first:
	\begin{equation}
		\begin{aligned}
		{[J_1}_\Lambda SB^\alpha] &= i (\chi + S) SB^\alpha, \\
		{[J_1}_\Lambda \omega_{\alpha \bar{\beta}}] &= - i {\omega_{\alpha
		\bar{\beta}}}_{,\gamma} SB^\gamma + i {\omega_{\alpha,
		\bar{\beta}}}_{,\bar{\gamma}} SB^{\bar{\gamma}} \\ 
		&= - i \Gamma^{\varepsilon}_{\alpha\gamma} \omega_{\varepsilon
		\bar{\beta}} SB^\gamma + i
		\Gamma^{\bar{\varepsilon}}_{\bar{\beta}\bar{\gamma}}
		\omega_{\alpha\bar{\varepsilon}} SB^{\bar{\gamma}},
		\end{aligned}
		\label{eq:prueba1}
	\end{equation}
	and with this we can compute easily:
	\begin{multline*}
		{[J_1}_\Lambda SB^\alpha SB^{\bar{\beta}}] = i \left((\chi + S)SB^\alpha
		\right) SB^{\bar{\beta}} + i SB^\alpha (\chi + S)
		SB^{\bar{\beta}} =\\ = i TB^\alpha SB^{\bar{\beta}} + i SB^\alpha
		TB^{\bar{\beta}}.
	\end{multline*}
	Therefore
	\begin{equation*}
		{[J_1}_\Lambda \Omega] = - i\Gamma^{\varepsilon}_{\alpha\gamma}
		\omega_{\varepsilon\bar{\beta}} SB^\gamma SB^\alpha
		SB^{\bar{\beta}} + i \omega_{\alpha\bar{\beta}} TB^\alpha
		SB^{\bar{\beta}} +  c.c.   = -g_{ij} TB^i SB^j,
	\end{equation*}
	and combining with (\ref{eq:conbeta}) we obtain:
	\begin{equation}
		{[J_1}_\Lambda J_2] = - \Gamma^k_{jl} g^{ij} SB^l(\Psi_i \Psi_k) -
		g^{ij} S\Psi_i \Psi_j - g_{ij}TB^i SB^j=: - H'.
		\label{eq:htildedef}
	\end{equation}
\end{proof}
\begin{rem}
	We note that there are no issues with quasi-associativity in the
	second term of (\ref{eq:htildedef}), that is why we did not include any parenthesis.

	Cubic expressions like this appeared in the literature linked to the
	Hamiltonian for the $N=2,2$ supersymmetric sigma model with target an
	untwisted generalized K\"ahler manifold. See for example \cite{bredthauer}
	and references therein. 
	\label{rem:citephysics}
\end{rem}
\begin{proof}[Proof of Theorem \ref{thm:n=2,2}] \hfill

	The main difficulty is trying to compute ${[H'}_\Lambda H']$ because of the
	term involving explicitly the Christoffel symbols. For this we will use the
	Jacobi identity for $N_K=1$ SUSY Lie Conformal algebras \cite{heluani3}. We
	have:
	\begin{equation*}
		\begin{aligned}
			{[ J_1}_\Lambda {[J_1}_\Gamma J_2]] &= - [
			[{J_1}_\Lambda J_1]_{\Lambda + \Gamma} J_2] -
			[{J_1}_\Gamma [{J_1}_\Lambda J_2]] \\ 
			- [{J_1}_\Lambda H'] &= [{H}_{\Lambda + \Gamma} J_2] +
			[{J_1}_\Gamma H'] \\ 
			- [{J_1}_\Lambda H'] - [{J_1}_\Gamma H'] &= (2T+ 2(\lambda + \gamma) + (\chi +
			\eta) S) J_2.
		\end{aligned}
	\end{equation*}
	From where we deduce:
	\begin{equation*}
		[{J_1}_\Lambda H'] = - (T + 2\lambda + \chi S) J_2,
	\end{equation*}
	therefore, using skew-symmetry:
	\begin{equation}
		[{H'}_\Lambda J_1] = (2T + 2\lambda + \chi S)J_2.
		\label{eq:nosecual3}
	\end{equation}
	Similarly, we find 
	\begin{equation}
		{[H'}_\Lambda J_2] = (2T + 2\lambda + \chi S)J_1.
		\label{eq:nosecual4}
	\end{equation}
	With these we can compute using the Jacobi identity again: 
	\begin{equation}
		\begin{aligned}
			{[H'}_\Lambda [{J_1}_\Gamma J_2]] &=  [ [{H'}_\Lambda
			J_1]_{\Lambda + \Gamma} J_2] + [{J_1}_\Gamma [{H'}_\Lambda
			J_2]] \\
			- [{H'}_\Lambda H'] &= [ (2T + 2\lambda + \chi S)
			{J_2}_{\Lambda + \Gamma} J_2] + [{J_1}_\Gamma (2 T + 2
			\lambda + \chi S) J_1]  \\ &= (- 2\gamma - \chi (\chi +
			\eta)) [{J_2}_{\Lambda+\Gamma}J_2] + \\ 
			& \quad +(2T + 2\gamma + 2
			\lambda + \chi (\eta + S)) [{J_1}_\Gamma J_1]   \\
			&= (2 \gamma -\lambda + \chi \eta) \left( H +
			2n (\chi + \eta)(\lambda + \gamma) \right) - \\ & \quad - (2 T +
			2\gamma + 2\lambda + \chi \eta + \chi S) (H + 2n \eta
			\gamma)  \\ [{H'}_\Lambda H'] &= (2T + 3\lambda + \chi S)
			H + 2n \lambda^2 \chi.
	\end{aligned}
	\label{eq:nosecual5}
	\end{equation}
	Similarly, we find
	\begin{equation}
		\begin{aligned}
			{[H}_\Lambda [{J_1}_\Gamma J_2] &= [(2 T + 2 \lambda + \chi
			S) {J_1}_{\Lambda + \Gamma} J_2] + {[J_1}_\Gamma (2 T+ 2
			\lambda + \chi S) J_2] \\ -{[H}_\Lambda H'] &= - (-2\gamma -
			\chi (\chi + \eta)) H' - (2T + 2 \gamma  + 2\lambda  + \chi
			(\eta + S)) H' \\
			{[H}_\Lambda H'] &= (2T + 3\lambda + \chi S) H'.
		\end{aligned}
		\label{eq:nosecual6}
	\end{equation}
	The theorem follows easily from equations
	(\ref{eq:nosecual3})-(\ref{eq:nosecual6}).
\end{proof}
\begin{proof}[Proof of Theorem \ref{thm:n=4,4}]\hfill \linebreak
		1. and 2.
	The fact that each pair $\{H, J_i\}$, $\{H, J_{\omega_i}\}$ generates an
	$N=2$ super vertex algebra of central charge $c = 12n$. Follows from
	\cite[Thm 7.4]{heluani2} and Theorem \ref{thm:anothern=2}. We need to  
	compute the commutation relations between the currents. Let us pick
	holomorphic coordinates for the first complex structure, so that in these
	coordinates, $J_1$ looks like:
	\begin{equation*}
		 J_1= i SB^\alpha \Psi_\alpha - i SB^{\bar\alpha} \Psi_{\bar\alpha} + i
         \textbf{g}_{,\alpha} TB^\alpha - i \textbf{g}_{,\bar\alpha}
         TB^{\bar\alpha}.
	\end{equation*}
	The other two K\"ahler forms combine to define a holomophic symplectic form
	$\eta = \omega_1 - i \omega_2$. It follows that the current $J^\pm :=
	\tfrac{1}{2} (J_{\omega_1} \mp i J_{\omega_2})$ is expressed in these
	coordinates as
	\begin{equation*}
		J^+ = \frac{1}{2} \left( \eta_{\alpha\beta} SB^\alpha SB^\beta +
		\eta^{\bar{\alpha}\bar{\beta}} \Psi_{\bar{\alpha}} \Psi_{\bar{\beta}}
		\right), \qquad J^- = \overline{J^+}.
	\end{equation*}
	We want to compute ${[J_1}_\Lambda J^\pm]$. For this we need
	(\ref{eq:prueba2}), (\ref{eq:prueba1}) and 
	\begin{equation*}
		{[J_1}_\Lambda \eta_{\alpha\beta}] = - i S(\eta_{\alpha\beta}),
		\qquad {[J_1}_\Lambda \eta^{\bar{\alpha}\bar{\beta}}] = i S
		(\eta^{\bar{\alpha}\bar{\beta}}),
	\end{equation*}
	to compute
	\begin{multline*}
		{[J_1}_\Lambda SB^\alpha SB^\beta ] = i \left( (\chi + S) SB^\alpha
		\right)SB^\beta - i SB^\alpha (\chi + S) SB^\beta = \\ = i (2\chi +
		S)
		SB^\alpha SB^\beta,
	\end{multline*}
	therefore 
	\begin{equation*}
		{[J_1}_\Lambda \eta_{\alpha \beta}SB^\alpha SB^\beta] = - i
		S(\eta_{\alpha\beta}) SB^\alpha SB^\beta + i \eta_{\alpha\beta}
		(2\chi + S) SB^\alpha SB^\beta,
	\end{equation*}
	and the first term vanishes since $\eta$ is closed, hence we obtain
	\begin{equation}
	 {[J_1}_\Lambda \eta_{\alpha \beta}SB^\alpha SB^\beta] = i (S + 2\chi)
	 \eta_{\alpha\beta} SB^\alpha SB^\beta. 	
	 \label{eq:prueban41}
	\end{equation}
	Similarly, we obtain
	\begin{equation*}
		\begin{aligned}
		{[J_1}_\Lambda \Psi_{\bar\alpha}\Psi_{\bar{\beta}}] &= i \left(
		(\chi + S) \Psi_{\bar{\alpha}}
		\right) \Psi_{\bar{\beta}} + i \lambda \mathbf{g}_{,\bar{\alpha}}
		\Psi_{\bar{\beta}} - i \Psi_{\bar{\alpha}} (\chi + S)
		\Psi_{\bar{\beta}} - i \Psi_{\bar{\alpha}} \lambda
		\mathbf{g}_{,\bar{\beta}} \\ &= i (S + 2\chi) \Psi_{\bar{\alpha}}
		\Psi_{\bar{\beta}} + i \lambda \left( \mathbf{g}_{,\bar{\alpha}}
		\Psi_{\bar{\beta}} - \mathbf{g}_{,\bar{\beta}} \Psi_{\bar{\alpha}}
		\right). 
	\end{aligned}
	\end{equation*}
	Therefore
	\begin{multline}
		{[J_1}_\Lambda \eta^{\bar{\alpha} \bar{\beta}} \Psi_{\bar{\alpha}}
		\Psi_{\bar{\beta}} ] = i (S + 2\chi)  \eta^{\bar{\alpha}
		\bar{\beta}} \Psi_{\bar{\alpha}}
		\Psi_{\bar{\beta}} + i \lambda \eta^{\bar{\alpha}\bar{\beta}} \left( \mathbf{g}_{,\bar{\alpha}}
		\Psi_{\bar{\beta}} - \mathbf{g}_{,\bar{\beta}} \Psi_{\bar{\alpha}}
		\right) + \\ + i \int_0^\Lambda {[S
		\eta^{\bar{\alpha}\bar{\beta}}}_\Gamma
		\Psi_{\bar{\alpha}}\Psi_{\bar{\beta}}] d\Gamma.
		\label{eq:prueban42}
	\end{multline}
	In order to compute the integral term, we need
	\begin{equation*}
		\begin{aligned}
		 {[S
		\eta^{\bar{\alpha}\bar{\beta}}}_\Lambda
		\Psi_{\bar{\alpha}}\Psi_{\bar{\beta}}] &= \chi \left(
		{\eta^{\bar{\alpha}\bar{\beta}}}_{,\bar{\alpha}} \Psi_{\bar{\beta}}
		- {\eta^{\bar{\alpha}\bar{\beta}}}_{,\bar{\beta}} \Psi_{\bar{\alpha}}
		\right), \\  &= - \chi \left(
		\Gamma^{\bar{\alpha}}_{\bar{\alpha}\bar{\varepsilon}}
		\eta^{\bar{\varepsilon}\bar{\beta}} \Psi_{\bar{\beta}}  +
		\Gamma^{\bar{\beta}}_{\bar{\alpha}\bar{\varepsilon}}
		\eta^{\bar{\alpha}\bar{\varepsilon}} \Psi_{\bar{\beta}} -
		\Gamma^{\bar{\alpha}}_{\bar{\beta}\bar{\varepsilon}}
		\eta^{\bar{\varepsilon}\bar{\beta}} \Psi_{\bar{\alpha}} -
		\Gamma^{\bar{\beta}}_{\bar{\beta}\bar{\varepsilon}}
		\eta^{\bar{\alpha}\bar{\varepsilon}} \Psi_{\bar{\alpha}}\right), \\
		&= - \chi \eta^{\bar\varepsilon \bar\beta} \left
		({\mathbf{g}}_{,\bar{\varepsilon}} \Psi_{\bar{\beta}}  -
		{\mathbf{g}}_{\bar{\beta}} \Psi_{\bar\varepsilon}\right).
	\end{aligned}
	\end{equation*}
	and from this equation and (\ref{eq:prueban42}) we obtain easily
	\begin{equation}
			{[J_1}_\Lambda \eta^{\bar{\alpha} \bar{\beta}} \Psi_{\bar{\alpha}}
		\Psi_{\bar{\beta}} ] = i (S + 2\chi)  \eta^{\bar{\alpha}
		\bar{\beta}} \Psi_{\bar{\alpha}}
		\Psi_{\bar{\beta}}.
		\label{eq:prueban43}
	\end{equation} 
	Combining with (\ref{eq:prueban41}) and their complex
	conjugates we obtain
	\begin{equation}
		{[J_1}_\Lambda J^\pm] = \pm i (S + 2 \chi) J^\pm.
		\label{eq:prueban44}
	\end{equation}
	It follows in the same way as in the proof of Theorem \ref{thm:anothern=2}
	that 
	\begin{equation}
		{[J^\pm}_\Lambda J^\pm] = 0,
		\label{eq:prueban45}
	\end{equation}
	therefore to finish the proof we need to compute ${[J^\pm}_\Lambda J^\mp]$.
	We proceed as in the proof of Theorem \ref{thm:anothern=2}. Let us define
	\begin{equation*}
		\beta = \eta^{\alpha\beta} \Psi_\alpha \Psi_\beta, \qquad \Omega =
		\eta_{\alpha \beta} SB^\alpha SB^\beta,
	\end{equation*}
	and their corresponding complex conjugates $\bar{\beta}, \bar\Omega$, so
	that we have
	\begin{equation*}
		J^+ = \frac{1}{2} \left( \Omega + \bar\beta \right), \qquad J^- =
		\frac{1}{2} \left( \bar{\Omega} + \beta \right).
	\end{equation*}
	Clearly we have 
	\begin{equation*}
		{[\beta}_\Lambda \bar{\beta}] = [\Omega_\Lambda \bar{\Omega}] = 0
	\end{equation*}
	We also have
	\begin{equation*}
		{[\Psi_\alpha}_\Lambda SB^\gamma SB^\delta] = \chi \left(
		{\delta_\alpha}^\gamma SB^\delta - {\delta_\alpha}^\delta SB^\gamma
		\right),
	\end{equation*}
	hence
	\begin{equation*}
		\begin{aligned}
		{[\Psi_\alpha}_\Lambda \Omega] &= (\eta_{\gamma,\delta})_{,\alpha}
		SB^\gamma SB^\delta + 2 \chi \eta_{\alpha\delta} SB^\delta,  \\ &= -
		\Gamma^{\varepsilon}_{\gamma\alpha} \eta_{\varepsilon \delta}
		SB^\gamma SB^\delta - \Gamma^\varepsilon_{\alpha \delta}
		\eta_{\gamma \varepsilon} SB^\gamma SB^\delta + 2 \chi \eta_{\alpha
		\delta} SB^\delta, \\ & =  - 2
		\Gamma^{\varepsilon}_{\gamma\alpha} \eta_{\varepsilon \delta}
		SB^\gamma SB^\delta + 2 \chi \eta_{\alpha
		\delta} SB^\delta, 
	\end{aligned}
	\end{equation*}
	and by skew-symmetry:
	\begin{equation*}
		\begin{aligned}
		{[\Omega}_\Lambda \Psi_\alpha ] &= -2\chi \eta_{\alpha\delta}
		SB^\delta - 2
		{\eta_{\alpha\delta}}_{,\gamma} SB^\gamma SB^\delta - 2
		\eta_{\alpha\delta} TB^\delta - 2 \Gamma^\varepsilon_{\gamma\alpha}
		\eta_{\varepsilon\delta} SB^\gamma SB^\delta,
		\\ &= - 2 \chi \eta_{\alpha\delta} SB^\delta + 2
		\Gamma^\varepsilon_{\delta\gamma} \eta_{\alpha \varepsilon}
		SB^\gamma SB^\delta - 2\eta_{\alpha\delta} TB^\delta, \\
		&= - 2 \eta_{\alpha\delta} (\chi + S) SB^\delta.
	\end{aligned}
	\end{equation*}
	Therefore
	\begin{equation*}
		\begin{aligned}
			{[\Omega}_\Lambda \Psi_\alpha \Psi_\beta] &= - 2
			\left( \eta_{\alpha\delta} (\chi + S)SB^\delta
			\right)\Psi_\beta + 2
			\Psi_\alpha \left( \eta_{\beta\delta} (\chi + S) SB^\delta
			\right) - \\ & \quad - 2 \int_0^\Lambda
			{[\eta_{\alpha\delta} (\chi + S)SB^\delta}_\Gamma
			\Psi_\beta] d\Gamma, \\ &= - 2
			\left( \eta_{\alpha\delta} (\chi + S)SB^\delta
			\right)\Psi_\beta + 2
			\left( \eta_{\beta\delta} (\chi + S) SB^\delta
			\right) \Psi_\alpha - \\  & \quad - 2 \int_0^\Lambda
			{[\eta_{\alpha\delta} (\chi + S)SB^\delta}_\Gamma
			\Psi_\beta] d\Gamma - 2 \chi \int_{-\nabla}^0
			{[\Psi_{\alpha}}_\Lambda \eta_{\beta \delta} 
		          SB^\delta] d \Lambda, \\ &=  - 2
			\left( \eta_{\alpha\delta} (\chi + S)SB^\delta
			\right)\Psi_\beta + 2
			\left( \eta_{\beta\delta} (\chi + S) SB^\delta
			\right) \Psi_\alpha - \\ & \quad - 2 \chi \lambda \eta_{\alpha
			\beta} - 2 \chi T (\eta_{\beta \alpha}),
		\end{aligned}
	\end{equation*}
	from where we deduce
	\begin{equation*}
		\begin{aligned}
			{[\Omega}_\Lambda \beta] &=  - 2 \eta^{\alpha \beta} \Bigl(
			\left( \eta_{\alpha\delta} (\chi + S)SB^\delta
			\right)\Psi_\beta \Bigr) + 2 \eta^{\alpha \beta} \Bigl(
			\left( \eta_{\beta\delta} (\chi + S) SB^\delta
			\right) \Psi_\alpha \Bigr) - \\ & \quad - 4n \chi \lambda
			- 2 \eta^{\alpha\beta} \chi T (\eta_{\beta \alpha}), \\ &=
			- 4 TB^\alpha \Psi_\alpha - 4 \chi SB^\alpha \Psi_\alpha +
			2 \chi T(\eta^{\alpha\beta}) \eta_{\alpha\beta} - 4 n \chi
			\lambda,
		\end{aligned}
	\end{equation*}
	and by skew-symmetry
	\begin{multline*}
		{[\beta}_\Lambda \Omega] =   + 4 \chi
		SB^\alpha \Psi_\alpha - 4 SB^\alpha S\Psi_\alpha  - 2 \chi
		T(\eta^{\alpha\beta}) \eta_{\alpha\beta} - \\ -2
		TS(\eta^{\alpha\beta})\eta_{\alpha\beta} - 2 T(\eta^{\alpha\beta})
		S(\eta_{\alpha\beta}) - 4n \chi \lambda.
	\end{multline*}
	A simple computation using that $\eta$ is parallel and that on a K\"ahler manifold we
	have $\Gamma_{\alpha\beta}^\beta = \mathbf{g}_{,\alpha}$ shows that
	\begin{equation*}
		T(\eta^{\alpha\beta}) \eta_{\alpha\beta} = - 2 \mathbf{g}_{,\gamma}
		TB^\gamma.
	\end{equation*}
	Hence collecting terms we obtain:
	\begin{equation}
		{[J^+}_\Lambda J^-] = - \frac{1}{2} \left( H +
		4n \lambda \chi \right) + \frac{\sqrt{-1}}{2} (S + 2\chi)
		J_1.
		\label{eq:prueban46}
	\end{equation}
	Equations (\ref{eq:prueban44}), (\ref{eq:prueban45}) and
	(\ref{eq:prueban46}) easily show that $\{H, J_1, J^+, J^-\}$ generate an
	$N=4$ super vertex algebra of central charge $c = 12 n$,
	thus proving 2).

       1) and 3). follow easily from 2). We can use the Jacobi identity for SUSY
       Lie conformal algebras to check 1) as follows. Let $H'_i = - {[J_i}_\Lambda
       J_{\omega_i}]$. By 2), $H'_1$ is half the $\eta$ coefficient in
       $-{[J_1}_\Lambda {[J_2}_\Gamma J_{\omega_3}]]$, but using the Jacobi
       identity, this is half the $\eta$ coefficient in 
       \begin{equation*}
	       {[ {[J_1}_\Lambda J_2]}_{\Lambda+\Gamma} J_{\omega_3}] +
	       {[J_2}_\Gamma {[J_1}_\Lambda J_{\omega_3}]].
       \end{equation*}
       Applying 2) again we can rewrite this as half the $\eta$ coefficient of
       \begin{equation*}
	       - (\eta - \chi) {[J_3}_{\Lambda + \Gamma} J_{\omega_3}] - (S + \eta
	       + 2\chi) {[J_2}_\Gamma J_{\omega_2}], 
       \end{equation*}
       which implies
       \begin{equation*}
	       H'_1 = \frac{1}{2} \left( H'_2 + H'_3 \right).
       \end{equation*}
       This equation together with its cyclic permutations imply $H'_1 = H'_2=H'_3$
       and therefore 1).

        To prove 3) we see that the fact that each $\{H^\pm, J_i^\pm\}$ are
	two commuting pairs of $N=2$ super vertex algebras of central charge $c = 6
	n$
	follows from Theorem \ref{thm:n=2,2}. To check that
	indeed we have two $N=4$ structures, we compute 
	\begin{multline*}
		{[J_1^\pm}_\Lambda J_2^\pm] = \frac{1}{4} {[J_1 \pm
		J_{\omega_1}}_\Lambda J_2 \pm J_{\omega_2}] = \\ =\frac{1}{4} \Bigl(
		[{J_1}_\Lambda J_2] + [{J_{\omega_1}}_\Lambda J_{\omega_2}] \pm
		[{J_1}_\Lambda J_{\omega_2}] \pm {[J_{\omega_1}}_\Lambda
		J_2]\Bigr) = \\ = \frac{1}{4} \Bigl( (S + 2\chi) J_3 + (S+2\chi)
		J_3 \pm (S + 2\chi) J_{\omega_3} \pm (S + 2\chi) J_{\omega_3}
		\Bigr) = \\ = (S + 2\chi) J_3^\pm.
	\end{multline*}
	Similarly, 
	\begin{multline*}
		{[J_1^\pm}_\Lambda J_2^\mp] = \frac{1}{4} {[J_1 \pm
		J_{\omega_1}}_\Lambda J_2 \mp J_{\omega_2}] = \\ =\frac{1}{4} \Bigl(
		[{J_1}_\Lambda J_2] - [{J_{\omega_1}}_\Lambda J_{\omega_2}] \mp
		[{J_1}_\Lambda J_{\omega_2}] \pm {[J_{\omega_1}}_\Lambda
		J_2]\Bigr) = \\ = \frac{1}{4} \Bigl( (S + 2\chi) J_3 - (S+2\chi)
		J_3 \mp \\ \mp  (S + 2\chi) J_{\omega_3} \pm (S + 2\chi) J_{\omega_3}
		\Bigr) = 0,
	\end{multline*}
	from where the Theorem follows.
\end{proof}


\begin{thebibliography}{10}

\bibitem{alekseev}
A.~Alekseev and T.~Strobl.
\newblock Current algebras and differential geometry.
\newblock {\em J. High Energy Phys.}, (3):035, 14 pp. (electronic), 2005.

\bibitem{heluani2}
D.~Ben-Zvi, R.~Heluani, and M.~Szcezny.
\newblock Supersymmetry of the chiral de {R}ham complex.
\newblock {\em Compositio Mathematicae}, 144(2):495--502, 2008.

\bibitem{bredthauer}
A.~Bredthauer, U.~Lindstr{\"o}m, J.~Persson, and M.~Zabzine.
\newblock Generalized {K}\"ahler geometry from supersymmetric sigma models.
\newblock {\em Lett. Math. Phys.}, 77(3):291--308, 2006.

\bibitem{bressler1}
P.~Bressler.
\newblock Vertex {A}lgebroids {I}.
\newblock {\em math/0202185}, 2002.

\bibitem{FL}
E.~Frenkel and A.~Losev.
\newblock Mirror symmetry in two steps: {A}-{I}-{B}.
\newblock {\em Comm. Math. Phys.}, 269(1):39--86, 2007.

\bibitem{frenkelnekrasov}
E.~Frenkel, A.~Losev, and N.~Nekrasov.
\newblock Instantons beyond topological theory {II}.
\newblock {\em preprint. math/0803.3302}, 2008.

\bibitem{gerbes2}
V.~Gorbounov, F.~Malikov, and V.~Schechtman.
\newblock Gerbes of chiral differential operators. {II}. {V}ertex algebroids.
\newblock {\em Invent. Math.}, 155(3):605--680, 2004.

\bibitem{gualtieri1}
M.~Gualtieri.
\newblock Generalized complex geometry.
\newblock {\em Oxford Ph.D. thesis}, 2004.

\bibitem{guttenberg}
S.~Guttenberg.
\newblock Brackets, sigma models and integrability of generalized complex
  structures.
\newblock {\em J. High Energy Phys.}, (6):004, 67 pp. (electronic), 2007.

\bibitem{heluani3}
R.~Heluani and V.G. Kac.
\newblock Super symmetric vertex algebras.
\newblock {\em Communications in mathematical physics}, (271):103--178, 2007.

\bibitem{kac:vertex}
V.~G. Kac.
\newblock {\em Vertex algebras for beginners}, volume~10 of {\em University
  Lecture}.
\newblock American Mathematical Society, 1996.

\bibitem{kacwakimoto1}
Victor~G. Kac and Minoru Wakimoto.
\newblock Quatum reduction and representation theory of superconformal
  algebras.
\newblock 2003.
\newblock math-ph/0304011.

\bibitem{kapustin}
A.~Kapustin.
\newblock Chiral de {R}ham complex and the half-twisted sigma model.
\newblock {\em preprint, hep-th/0504074}, 2005.

\bibitem{malikov}
A.~Malikov, V.~Shechtman, and A.~Vaintrob.
\newblock Chiral de {R}ham complex.
\newblock {\em Comm. Math. Phys}, 204(2):439--473, 1999.

\bibitem{Witten}
E.~Witten.
\newblock Two-dimensional sigma models with $(0,2)$ supersymmetry:{P}ertubative
  aspects.
\newblock {\em preprint, hep-th/0504078}, 2005.

\end{thebibliography}
\end{document}